\documentclass[a4paper, 11pt, final]{article}

\usepackage{enumerate}
\usepackage{amsfonts}
\usepackage{graphicx}
\usepackage{psfrag}
\usepackage{array}

\usepackage{hyperref}

\newtheorem{proposition}{Proposition}[section]
\newtheorem{theorem}[proposition]{Theorem}
\newtheorem{lemma}[proposition]{Lemma}
\newtheorem{corollary}[proposition]{Corollary}
\newtheorem{definition}[proposition]{Definition}

{\hspace{-5pt}{\nobreak\quad\nobreak\hfill\nobreak$\square$\vspace{8pt}%
\par}\smallskip\goodbreak}

\newenvironment{proofof}[1]{\smallskip\noindent\emph{Proof of #1.}%
\hspace{1pt}}{\hspace{-5pt}{\nobreak\quad\nobreak\hfill\nobreak%
$\square$\vspace{8pt}\par}\smallskip\goodbreak}

\newcommand{\Section}[1]{\section{#1}\setcounter{equation}{0}}

\renewcommand{\phi}{\varphi}
\renewcommand{\L}[1]{\mathbf{L^#1}}

\newcommand{\Lloc}[1]{\mathbf{L^{#1}_{loc}}}
\newcommand{\C}[1]{\mathbf{C^{#1}}}
\newcommand{\Cc}[1]{\mathbf{C_c^{#1}}}

\newcommand{\modulo}[1]{{\left|#1\right|}}
\newcommand{\norma}[1]{{\left\|#1\right\|}}

\newcommand{\reali}{{\mathbb{R}}}

\newcommand{\tv}{\mathrm{TV}}
\newcommand{\BV}{\mathbf{BV}}

\renewcommand{\epsilon}{\varepsilon}
\newcommand{\pint}[1]{\mathaccent23{#1}}

\newcommand{\fu}[1]{\hat\phi_{#1}}
\newcommand{\fd}[1]{\check\phi_{#1}}

\title{On the 1D Modeling \\ of Fluid Flowing through a Junction}

\author{Rinaldo M.~Colombo \\ Dipartimento di Matematica \\ 25133
  Brescia \\ Italy \and Mauro Garavello\thanks{Partially supported
    by Dipartimento di Matematica e Applicazioni of the University of
    Milano-Bicocca} \\ Di.S.T.A. \\ 15100
  Alessandria \\ Italy}

\begin{document}

\maketitle

\begin{abstract}
  Consider a fluid flowing through a junction between two pipes with
  different sections. Its evolution is described by the 2D or 3D Euler
  equations, whose analytical theory is far from complete and whose
  numerical treatment may be rather costly. This note compares
  different 1D approaches to this phenomenon.

  \medskip

  \noindent\textit{2000~Mathematics Subject Classification:} 35L65
  
  \medskip
  
\end{abstract}

\Section{Introduction}
\label{sec:Intro}

The motion of a fluid in a pipe is described by the one dimensional
$p$-system in Eulerian coordinates, which reads
\begin{equation}
  \label{eq:System}
  \left\{
    \begin{array}{lcl}
      \displaystyle
      \partial_t \rho + \partial_x q & = & 0
      \\
      \displaystyle
      \partial_t q + \partial_x \left( \frac{q^2}{\rho} + p(\rho) \right)
      & = & 0 \,,
    \end{array}
  \right.
\end{equation}
where $\rho$ is the fluid mass density, $q$ its linear momentum
density, $p$ its pressure, $x$ is the space coordinate and $t$ is
time.  Assume that the section of the pipe has a change due to a
junction, sited at, say, $x=0$. The motion of the fluid can then be
described by~(\ref{eq:System}) together with a condition on the traces
of the thermodynamic variables at $x=0$, i.e.~a condition of the type
\begin{equation}
  \label{eq:Phi}
  \Phi 
  \left( 
    a_l,\rho(t, 0-), q(t, 0-); a_r, \rho(t, 0+), q(t,0+)
  \right) 
  = 0 \,;
\end{equation}
here $a_l$ and $a_r$ are the pipes' sections to the left and to the
right of the junction.  Various choices of the function $\Phi$ are
present in the literature, see for instance~\cite{BandaHertyKlar2,
  BandaHertyKlar1, ColomboGaravello1, ColomboGaravello2,
  ColomboMarcellini, GuerraMarcelliniSchleper}. Here, we compare the
various definitions on the basis, \emph{a priori}, of their physical
meaning and, \emph{a posteriori}, of their analytical properties.

This phenomenon is essentially 3D, due to turbulence, or may be
reduced to 2D, but it is intrinsically multi-dimensional: \emph{``It
  is probable that, in the future, changes of cross-sectional area in
  pipework systems will be incorporated as multi-dimensional patches
  in one-dimensional simulations.''}, \cite[Section~6.1.,
p.~262]{WinterbonePearson}. However, the numerical treatment of 2D or
3D hyperbolic systems is far more expensive than that of 1D
systems. We refer to~\cite[Table~1]{Herty} for a striking comparison
between CPU times in 2D and 1D numerical integrations of Euler
equations in the case of a ``T'' junction. Besides, as is well known,
basic analytical questions, such as well posedness, about 2D or 3D
Euler equations are still unanswered.  Aiming also at the case of more
complex gas networks, the availability of 1D simple, though
approximate, \emph{``patches''} (with reference to the citation above)
may be of great help. All this implies, in turn, to choose a specific
function $\Phi$ in~(\ref{eq:Phi}). Below, we compare various choices
appeared in the recent literature.

First, we deal with solutions to Riemann Problems
for~(\ref{eq:System}) at a junction between two pipes with different
sections. By Riemann Problem for~(\ref{eq:System}) we mean
\begin{equation}
  \label{eq:RP}
  \left\{
    \begin{array}{l}
      \partial_t \rho + \partial_x q = 0,
      \\
      \partial_t q + \partial_x \left(q^2 /\rho + p(\rho)\right) = 0,
      \\
      (\rho,q) (0,x) = 
      \left\{
        \begin{array}{ll}
          (\bar \rho_l,\bar q_l), & \textrm{ if } x < 0,
          \\
          (\bar \rho_r,\bar q_r), & \textrm{ if } x > 0,
        \end{array}
      \right.
    \end{array}
  \right.
\end{equation}
where $(\bar \rho_l,\bar q_l)$ and $(\bar \rho_l,\bar q_l)$ are fixed
in $\pint{\reali}^+ \times \reali$. Recall that the Riemann Problem
plays an essential role in the construction of solutions to the
general Cauchy Problem
for~(\ref{eq:System})--(\ref{eq:Phi}). Analogously, the solution
of~(\ref{eq:System})--(\ref{eq:Phi}) at a junction is essential for
the development of the theory of the $p$-system on general networks,
see~\cite{GaravelloPiccoli}.

Following various results in the literature, such
as~\cite{BandaHertyKlar2, BandaHertyKlar1, ColomboGaravello1,
  ColomboGaravello2, ColomboGaravello3, ColomboHertySachers,
  ColomboMarcellini, HertySeaid, HoldenRisebroKink}, we restrict our
attention to subsonic flows in both pipes.

Remark that the modeling of compressible fluid flowing through pipes
with sudden enlargements can be achieved through several entirely
different techniques, for instance the isothermal lattice-Boltzmann
methods, see~\cite{Agrawal}, or the full system of Navier--Stokes
equations, see~\cite{Duda}, whereas experimental data are also
available in the literature, see~\cite{Dekam, cinesi, Rathakrishnan}.

The next section is devoted to~(\ref{eq:RP}) and to the Cauchy problem
for~(\ref{eq:System})--(\ref{eq:Phi}), both in general and for
specific choices of the junction condition $\Phi$. Then,
Section~\ref{sec:Numerical} displays some numerical solutions. All the
analytical details are gathered in Section~\ref{sec:Tech}.

%
%
%
%
\Section{Analytical Results}
\label{sec:Nota}

Throughout, we identify $u \equiv (\rho,q)$. On the pressure law $p$,
we assume
\begin{description}
\item[(EoS)] $p \in \C2(\reali^+; \reali^+)$ is such that for all
  $\rho > 0$, $p'(\rho) >0$ and $p''(\rho) \geq 0$.
\end{description}
\noindent The classical example is the $\gamma$-law, where $p(\rho) =
k \, (\rho/\bar\rho)^\gamma$, where $\gamma \geq 1$ and $\bar\rho >
0$. Recall for later use the expression of the flow of the linear
momentum:
\begin{displaymath}
  P(u) 
  =
  \frac{q^2}\rho+p(\rho) \,.
\end{displaymath}
Throughout, the present analysis is limited to the \textit{subsonic}
region given by
\begin{equation}
  \label{subsonic2}
  A_{0} 
  =
  \left\{ 
    u \in \pint{\reali}^{+} \times \reali
    \colon\quad
    \lambda_{1} (u) < 0 <\lambda_{2} (u) 
  \right\} \,,  
\end{equation}
where $\lambda_i$ is the $i$-th eigenvalue of~(\ref{eq:System}),
see~(\ref{eigenvectors}).  Recall the following definition of solution
to the Riemann Problem~(\ref{eq:RP}) under condition~(\ref{eq:Phi}),
see also~\cite[Definition~2.2]{ColomboMarcellini}.

%
%
\begin{definition}
  \label{def:PhiSol}
  A weak $\Phi$-solution to the Riemann Problem~(\ref{eq:RP}) is a map
  \begin{equation}
    \label{eq:regularity}
    \begin{array}{rcl}
      u
      & \in & 
      \C0 \left( 
        \reali^+; 
        \Lloc1 (\reali; \pint{\reali}^+ \times \reali) \right)
      \\ 
      u (t) & \in & \BV (\reali; \pint{\reali}^+ \times \reali)
      \quad \mbox{ for a.e. } t \in \reali^+
    \end{array}
  \end{equation}
  such that
  \begin{enumerate}

  \item for all $\phi \in \Cc1 (\pint{\reali}^+ \times \reali;
    \reali)$ whose support does not intersect $x=0$
    \begin{displaymath}
      \int_{\reali^+} \int_{\reali} 
      \left( 
        \left[
          \begin{array}{c}
            \rho \\ q
          \end{array}
        \right]
        \, \partial_t \phi + 
        \left[ 
          \begin{array}{c}
            q
            \\
            P(\rho,q)
          \end{array}
        \right]
        \, \partial_x \phi
      \right)
      \, dx \, dt = 0 \,;
    \end{displaymath}

  \item for a.e.~$x \in \reali$, $u(0,x)$ coincides with the initial
    datum in~(\ref{eq:RP});
 
  \item for a.e.~$t \in \reali^+$, the junction
    condition~(\ref{eq:Phi}) at the junction is met.
  \end{enumerate}
\end{definition}

\noindent We consider the following properties of the junction
condition~(\ref{eq:Phi}), which we rewrite here as
\begin{equation}
  \label{eq:Sigma}
  \Phi(a_l, u_l; a_r, u_r) =0 \,.
\end{equation}
\begin{enumerate}[\bf($\mathbf{\Phi}$1)]\setcounter{enumi}{-1}
\item \label{it:p0} Regularity: $\Phi \in \C1 \left( (\pint{\reali}^+
    \times A_0 )^2; \reali^2 \right)$, where $A_0$ is given
  by~(\ref{subsonic2}). Moreover, the $2\times 2$ matrix $D_{u_r}
  \Phi(a_l,u_l;a_r,u_r)$ is invertible, for all $a_l, a_r >0$ and
  $u_l,u_r \in \pint{A}_0$.

\item \label{it:p1} No-junction case: for all $a>0$ and $u_l,u_r \in
  A_0$, $\Phi(a, u_l, a, u_r) = 0$ if and only if $u_l = u_r$.

\item \label{it:p2} Left-right symmetry: for all $a_l,a_r>0$ and
  $(\rho_l, q_l), (\rho_r, q_r) \in A_0$, $\Phi(a_l,\rho_l, q_l;
  a_r,\rho_r,q_r) = 0$ if and only if $\Phi(a_r, \rho_r, -q_r;a_l ,
  \rho_l, -q_l) = 0$.

\item \label{it:p3} Consistency: for all positive $a_l, a_m, a_r$ and
  all $u_l, u_m, u_r \in A_0$, if $\Phi(a_l,u_l;a_m,u_m)=0$ and
  $\Phi(a_m,u_m;a_r,u_r) = 0$ then $\Phi(a_l,u_l;a_r,u_r) = 0$.

\item \label{it:p4} Hydrostatic limit: for all positive $a_l, a_r$ and
  for all densities $\rho_l, \rho_r$, $\Phi(a_l,\rho_l,0;a_r,\rho_r,0)
  = 0$ if and only if $p(\rho_l) = p(\rho_r)$.
\end{enumerate}

\noindent Moreover, by an immediate extension
of~\cite[Lemma~2.1]{ColomboMarcellini},
\textbf{($\mathbf{\Phi}$\ref{it:p0})} ensures that~(\ref{eq:Sigma})
implicitly defines a map $u_r = T(u_l; a_l, a_r)$ in a neighborhood of
a subsonic state satisfying $\Phi(a_l,u_l;a_r,u_r)=0$.  In turn, this
implies the local well posedness of Cauchy problems for data near to
stationary solutions; see~\cite{ColomboGaravello3,
  ColomboHertySachers, ColomboMarcellini} and
Proposition~\ref{prop:Local} below.  Here, \emph{``near''} is meant in
the sense of the total variation.

\begin{proposition}
  \label{prop:Local}
  Let $p$ satisfy~\textbf{(EoS)}, $\Phi$
  satisfy~\textbf{($\mathbf{\Phi}$\ref{it:p0})}.  Then, for all
  positive $\bar a_l, \bar a_r$ and $\bar u_l, \bar u_r \in
  \pint{A}_0$ such that $\Phi(\bar a_l,\bar u_l; \bar a_r, \bar u_r) =
  0$, there exist neighborhoods $\mathcal{A}_l, \mathcal{A}_r$ of
  $\bar a_l, \bar a_r$ and $\mathcal{U}_l, \mathcal{U}_r$ of $\bar
  u_l, \bar u_r$ such that for all $a_l \in \mathcal{A}_l$, $a_r \in
  \mathcal{A}_r$ and $u_l \in \mathcal{U}_l$, there exists a unique
  $u_r \in \mathcal{U}_r$ such that the map
  \begin{displaymath}
    \hat u 
    = 
    u_l \, \chi_{\strut \left]-\infty, 0\right[} +
    u_r \, \chi_{\strut \left]0, +\infty\right[}
  \end{displaymath}
  is a stationary weak $\Phi$-solution
  to~(\ref{eq:System})--(\ref{eq:Sigma}) in the sense of
  Definition~\ref{def:PhiSol}. Moreover, there exist $\delta, L > 0$
  and a semigroup $S \colon \pint{\reali}^+ \times \mathcal{D} \to
  \mathcal{D}$ such that
  \begin{enumerate}
  \item $\mathcal{D} \supseteq \left\{ u \in \hat u + \L1(\reali; A_0)
      \colon \tv(u - \hat u) < \delta \right\}$.
  \item For all $u \in \mathcal{D}$, $S_0 u = u$ and for all $t,s \geq
    0$, $S_t S_s u = S_{s+t} u$.
  \item For all $u, \tilde u \in \mathcal{D}$ and for all $t,\tilde t
    \geq 0$,
    \begin{displaymath}
      \norma{S_tu - S_{\tilde t} \tilde u}_{\L1}
      \leq
      L \cdot \left(
        \norma{u - u}_{\L1} + \modulo{t-\tilde t}
      \right) .
    \end{displaymath}
  \item If $u \in \mathcal{D}$ is piecewise constant, then for $t$
    small, $S_t u$ is the gluing of Lax solutions to standard Riemann
    problems at the points of jump $x \neq 0$ in $u$, and of solutions
    to~(\ref{eq:RP})--(\ref{eq:Phi}) in the sense of
    Definition~\ref{def:PhiSol} at the junction at $x = 0$.
  \item For all $u \in \mathcal{D}$, the orbit $t \mapsto S_t u$ is a
    weak $\Phi$-solution to~(\ref{eq:System}).
  \item For any sequences $a_l^n \in \mathcal{A}_l$, $a_r^n \in
    \mathcal{A}_r$ and $u_l^n \in \mathcal{U}_l$, call $S^n$ the
    corresponding semigroup. If $a^n_l \to a_l \in \mathcal{A}_l$,
    $a^n_r \to a_r \in \mathcal{A}_r$ and $u^n_r \to u_l \in
    \mathcal{U}_l$, then $S^n$ converges uniformly on any compact time
    interval to the semigroup $S$ defined by $a_l, a_r$ and $u_l$ in
    $\Lloc1(\reali;A_0)$.
  \end{enumerate}
\end{proposition}

\noindent Suitable choices of $\Phi$ in~(\ref{eq:Sigma}) allow to
recover various definitions in the current literature.  In the
examples of the next table, the first component $\Phi_1$ of $\Phi$ is
equal to $a_l q_l - a_r q_r$, ensuring the conservation of mass.

\smallskip

\begin{center}
  \begin{tabular}{|m{0.05\linewidth}|m{0.34\linewidth}|m{0.48\linewidth}|}
    \hline
    &
    $\Phi_2(a_l,u_l,a_r,u_r)$
    &
    Meaning
    \\
    \hline
    \textbf{(L)}
    &
    $a_r P(u_r) - a_l P(u_l)$
    &
    Conservation of linear momentum, see~\cite{ColomboGaravello2}
    \\ 
    \cline{1-3}
    \textbf{(p)}
    &
    $p(\rho_r) - p(\rho_l)$
    &
    Equal pressure, typically motivated by static equilibrium,
    see~\cite{BandaHertyKlar2, BandaHertyKlar1}
    \\
    \hline
    \textbf{(P)}
    &
    $P(u_r) - P(u_l)$
    & 
    Equal dynamic pressure, see~\cite{ColomboGaravello1,
      ColomboGaravello3}
    \\
    \hline
    \textbf{(S)}
    &
    $
    \begin{array}{l}
      \displaystyle
      \!\!
      a_r P(u_r) - a_l P(u_l)\\
      -
      \int_{a_l}^{a_r} p \left( R(\alpha;\rho_l,q_l) \right) \, d\alpha\!\!
    \end{array}
    $
    &
    Limit of the condition for smooth variations of the pipes' sections,
    see~\cite{ColomboMarcellini, GuerraMarcelliniSchleper}
    \\
    \hline
  \end{tabular}
\end{center}
\medskip

\noindent The leftmost column displays the letter used below to refer
to the solutions of~(\ref{eq:System})--(\ref{eq:Phi}) with
$\Phi_1(a_l,\rho_l,q_l;a_r,\rho_r,q_r) = a_r q_r - a_l q_l$ and
$\Phi_2$ as in the second column. The function $R$ in~\textbf{(S)} is
defined in~(\ref{eq:sysR}), see also~(\ref{eq:SigmaAlfa}). We use the
above conditions to select particular weak solutions
to~(\ref{eq:System})--(\ref{eq:Phi}).

The main analytical properties of these choices are summarized below.

\begin{proposition}
  \label{prop:allP}
  All junction conditions $\Phi \colon ( {\reali}^+ \times A_0 )^2 \to
  \reali^2 $ such that $\Phi_1(a_l,\rho_l,q_l;a_r,\rho_r,q_r) = a_r
  q_r - a_l q_l$ and $\Phi_2$ is given by~\textbf{(L)}, \textbf{(p)},
  \textbf{(P)} or~\textbf{(S)} satisfy
  \textbf{($\mathbf{\Phi}$\ref{it:p0})},
  \textbf{($\mathbf{\Phi}$\ref{it:p1})},
  \textbf{($\mathbf{\Phi}$\ref{it:p2})} and
  \textbf{($\mathbf{\Phi}$\ref{it:p3})}.
  Property~\textbf{($\mathbf{\Phi}$\ref{it:p4})} is satisfied
  by~\textbf{(p)}, \textbf{(P)} and~\textbf{(S)}, not by~\textbf{(L)}.
\end{proposition}

Besides, below we characterize the set of initial data $\bar u_l$,
$\bar u_r$ such that~(\ref{eq:RP}) admits a solution and such that
this solution is unique.

Throughout, we refer to the forward, respectively, backward, Lax
curves of the $i$-th family $q = L_i(\rho;\rho_o,q_o)$, respectively
$q = L_i^-(\rho;\rho_o,q_o)$; see~(\ref{eq:LaxCurves}) for the
explicit expressions. As usual, $\lambda_i$ and $r_i$ are the $i$-th
eigenvalues and right eigenvectors. The explicit expressions are
collected in Section~\ref{sec:Tech}, together with several technical
proofs.
\begin{figure}[htpb]
  \centering
  \begin{psfrags}
    \psfrag{bar}{{\scriptsize $(\bar u)$}} \psfrag{q}{$q$}
    \psfrag{rho}{$\rho$} \psfrag{ul}{{\scriptsize $\fu{l} (\bar u)$}}
    \psfrag{ur}{{\scriptsize $\fu{r} (\bar u)$}}
    \psfrag{dl}{{\scriptsize $\fd{l} (\bar u)$}}
    \psfrag{dr}{{\scriptsize $\fd{r} (\bar u)$}}
    \includegraphics[width=6cm]{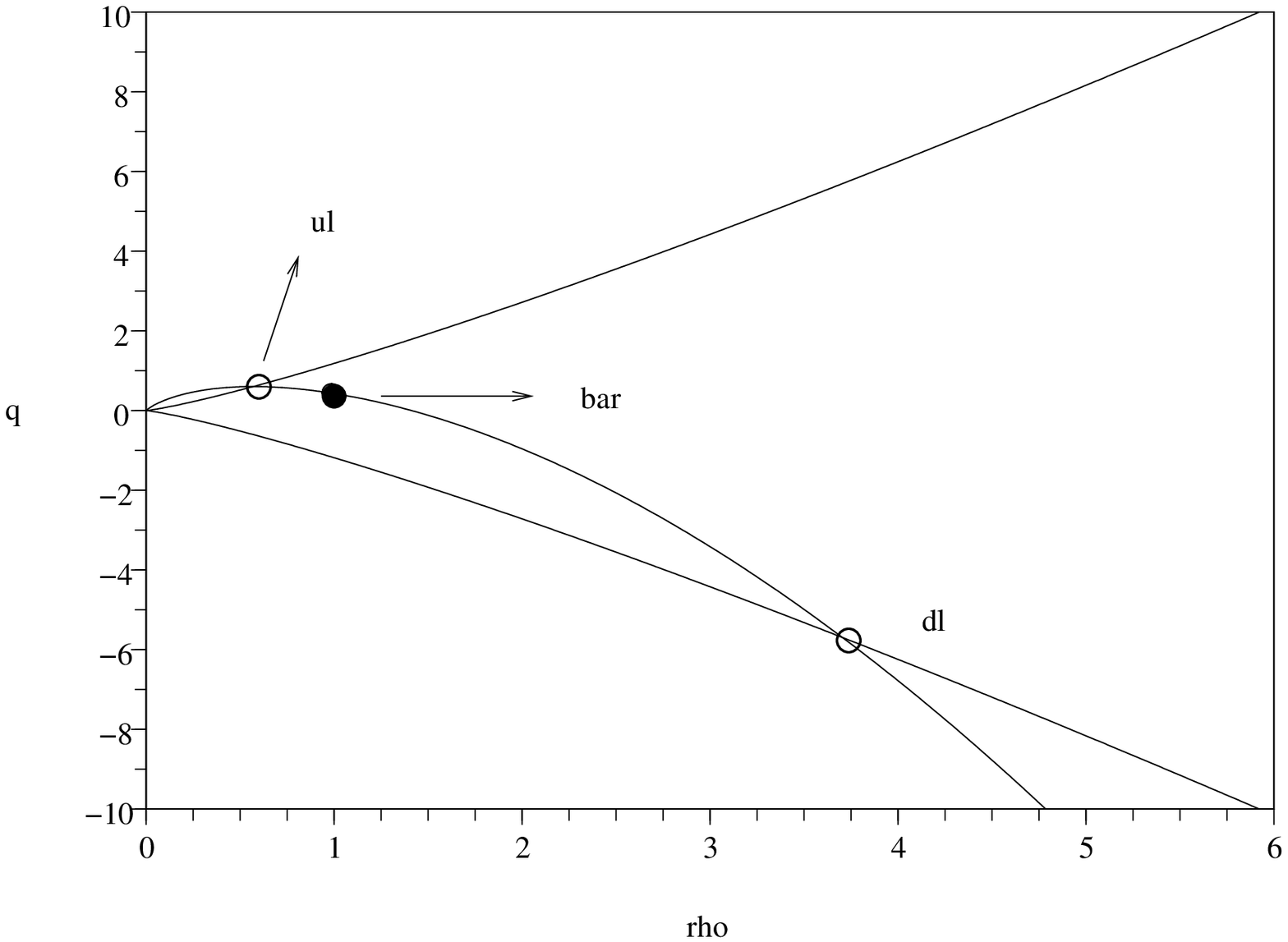}
  \end{psfrags}
  \begin{psfrags}
    \psfrag{bar}{{\scriptsize $(\bar u)$}} \psfrag{q}{$q$}
    \psfrag{rho}{$\rho$} \psfrag{ul}{{\scriptsize $\fu{l} (\bar u)$}}
    \psfrag{ur}{{\scriptsize $\fu{r} (\bar u)$}}
    \psfrag{dl}{{\scriptsize $\fd{l} (\bar u)$}}
    \psfrag{dr}{{\scriptsize $\fd{r} (\bar u)$}}
    \includegraphics[width=6cm]{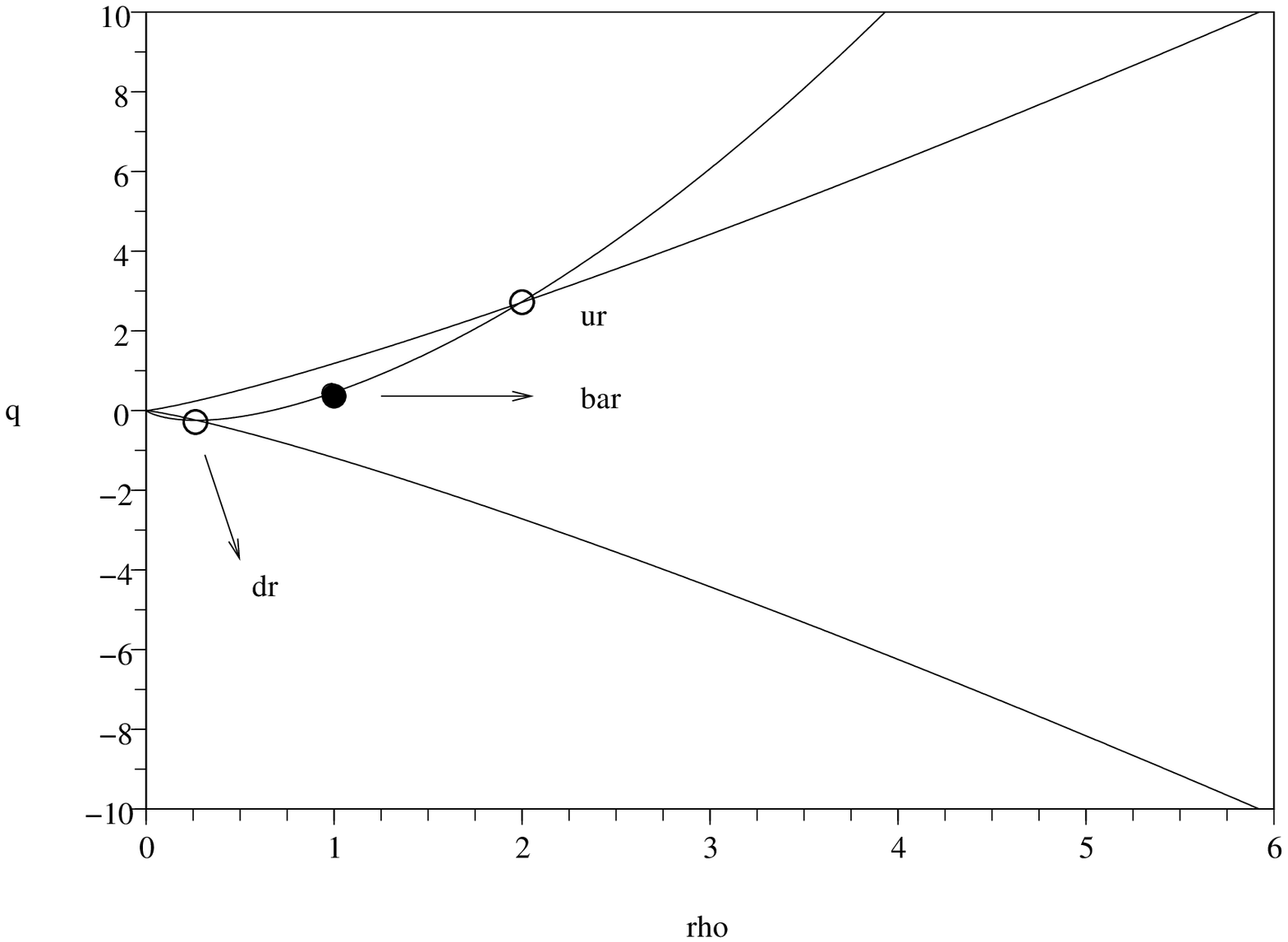}
  \end{psfrags}
  \caption{\small{Left, the densities $\fd{l} (\bar u)$ and $\fu{l}
      (\bar u)$, along a $1$--Lax curve; right, the densities $\fd{r}
      (\bar u)$ and $\fu{r} (\bar u)$ along a reversed $2$-Lax curve;
      all defined by~(\ref{eq:phi_eq}).}}
  \label{fig:fhi_hat}
\end{figure}
Moreover, for $\bar u \in A_0$, the densities $\fd{l} (\bar u)$,
$\fu{l} (\bar u)$, $\fd{r} (\bar u)$, $\fu{r} (\bar u)$ denote the
$\rho$-component of the intersections of the Lax curves $\rho \mapsto
\left( \rho, L_2^-(\rho; \bar u) \right)$ and $\rho \mapsto \left(
  \rho, L_1(\rho; \bar u) \right)$ with the boundary of the subsonic
region $A_0$. They are uniquely determined by:
\begin{equation}
  \label{eq:phi_eq}
  \!\!\!\!
  \begin{array}{rcl@{\quad}rcl}
    \lambda_2 \left( \fd{r} (\bar u) ,
      L_2^-\left( \fd{r} (\bar u); \bar u \right)\right)
    & = & 0,
    &
    \lambda_2 \left( \fd{l} (\bar u) ,
      L_1\left( \fd{l} (\bar u); \bar u \right) \right)
    & = & 0,
    \\
    \lambda_1 \left( \fu{r} (\bar u) ,
      L_2^-\left( \fu{r} (\bar u); \bar u \right) \right)
    & = & 0,
    &
    \lambda_1 \left( \fu{l} (\bar u) ,
      L_1\left( \fu{l} (\bar u); \bar u \right) \right)
    & = & 0.
  \end{array}
\end{equation}
see Figure~\ref{fig:fhi_hat}. Note that $\fu{l} (\bar u) \leq \fu{r}
(\bar u)$, $\fd{r} (\bar u) \leq \fd{l} (\bar u)$, $\fd{r} (\bar u)
\leq \fu{r} (\bar u)$, $\fu{l} (\bar u) \leq \fd{l} (\bar u)$.

%
%
%
%
\subsection{\textbf{(L)}-solutions}
\label{sec:L-sol}

Under condition~\textbf{(L)}, (\ref{eq:Sigma}) reads
\begin{displaymath}
  a_l \, q_l = a_r \, q_r
  \qquad\mbox{ and }\qquad
  a_l \, P(\rho_l, q_l)
  = 
  a_r \, P(\rho_r, q_r)
\end{displaymath}
and yields the standard Riemann solver giving Lax solutions. This is
an immediate extension of the standard Lax solutions
of~(\ref{eq:System}).

\begin{theorem}
  \label{thm:Lsol}
  Let~\textbf{(EoS)} hold and fix positive $a_l, a_r$. Consider the
  Riemann Problem~(\ref{eq:RP}) with $u_l, u_r \in A_0$. Define
  \begin{eqnarray*}
    l'
    & = &
    \inf \left\{ \rho \in [\fu{l}(\bar u_l),
      \fd{l}(\bar u_l)] \colon
      a_l \, L_1(\rho; \bar u_l) \le
      a_r \, L_2^-(\fu{r}(\bar u_r); \bar u_r)
    \right\}
    \\
    l''
    & = &
    \sup \left\{ \rho \in [\fu{l}(\bar u_l),
      \fd{l}(\bar u_l)] \colon
      a_l \, L_1(\rho; \bar u_l) \ge
      a_r \, L_2^-(\fd{r}(\bar u_r); \bar u_r)
    \right\}.
  \end{eqnarray*}
  A necessary and sufficient condition for the existence and
  uniqueness of an~\textbf{(L)}-solution to~(\ref{eq:RP}), attaining
  values in $A_0$, is
  \begin{equation}
    \label{eq:L-sol_CNS}
    \left\{
      \begin{array}{rcl}
        a_l P \left(l', L_1(l'; \bar u_l)\right) 
        & \leq &
        a_r P\left (g \left(\frac{a_l}{a_r}
            L_1(l'; \bar u_l) \right),
          \frac{a_l}{a_r} L_1(l'; \bar u_l)\right) 
        \\
        a_l P \left(l'', L_1(l''; \bar u_l)\right) 
        & \geq &
        a_ rP
        \left( g \left(\frac{a_l}{a_r} L_1(l''; \bar u_l)\right),
          \frac{a_l}{a_r} L_1(l''; \bar u_l)\right) 
      \end{array}
    \right.
  \end{equation}
  where $g$ is the inverse in $A_0$ of the map $\rho \to L_2^-(\rho ;
  \bar u_r)$.
\end{theorem}

\noindent Note that $l' \le l''$, since the sets defining $l'$ and
$l''$ are not disjoint.

An example of non existence of solutions to~(\ref{eq:RP}) is provided
by the following corollary.

\begin{figure}[htpb]
  \centering
  \begin{minipage}{0.45\linewidth}
    \begin{psfrags}
      \psfrag{rho}{$\rho$} \psfrag{q}{$q$} \psfrag{l}{{\tiny $(\bar
          u_l)$}} \psfrag{r}{{\tiny $(\bar u_r)$}} \psfrag{fl}{{\tiny
          $\fu{l}(\bar u_l)$}} \psfrag{fr}{{\tiny $\fu{r}(\bar u_r)$}}
      \includegraphics[width=5cm]{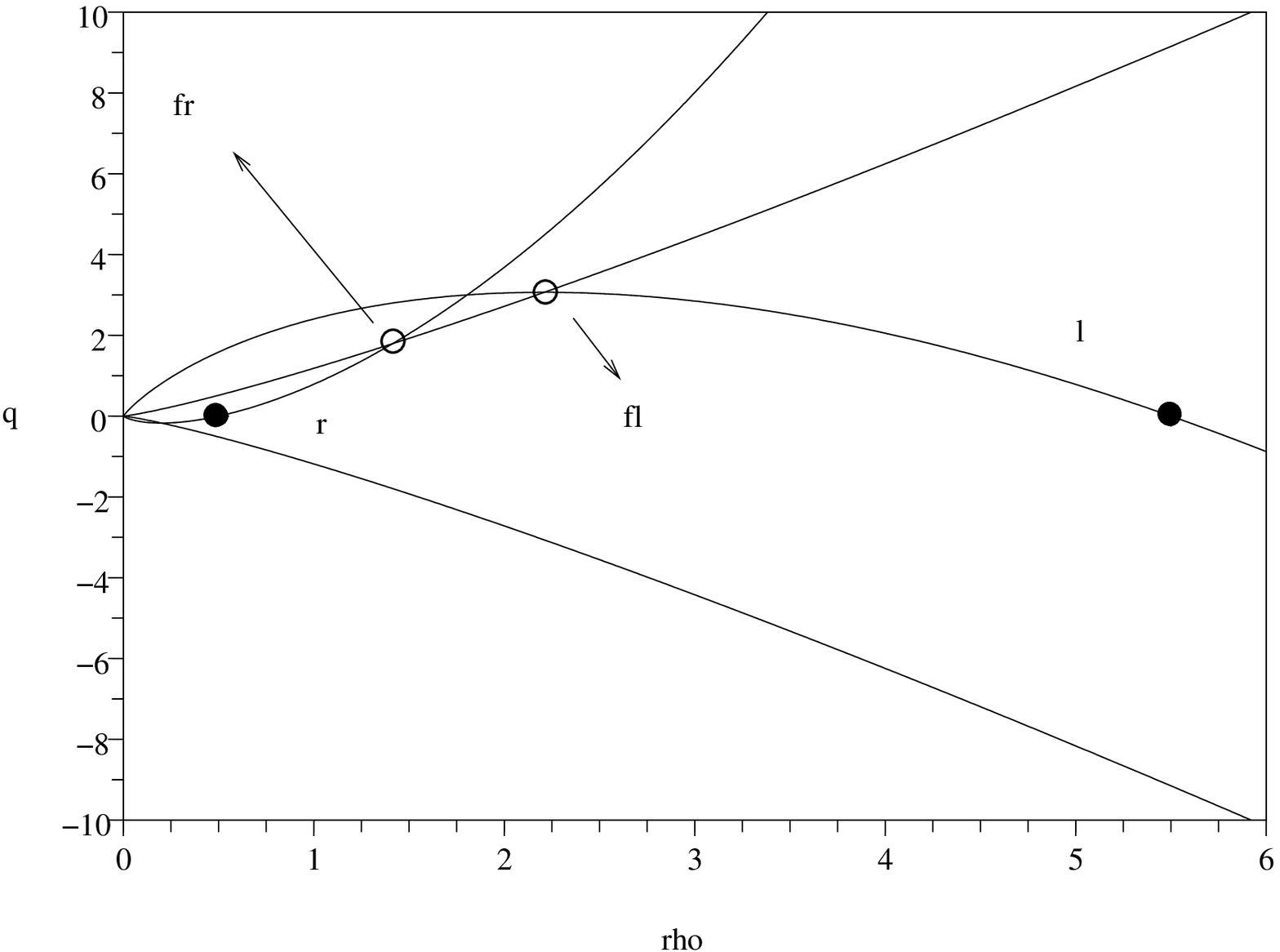}
    \end{psfrags}
  \end{minipage}
  \begin{minipage}{0.45\linewidth}
    \begin{psfrags}
      \psfrag{rho}{$\rho$} \psfrag{q}{$q$} \psfrag{l}{{\tiny $(\bar
          u_l)$}} \psfrag{r}{{\tiny $(\bar u_r)$}} \psfrag{fl}{{\tiny
          $\fu{l}(\bar u_l)$}} \psfrag{fr}{{\tiny $\fu{r}(\bar u_r)$}}
      \includegraphics[width=5cm]{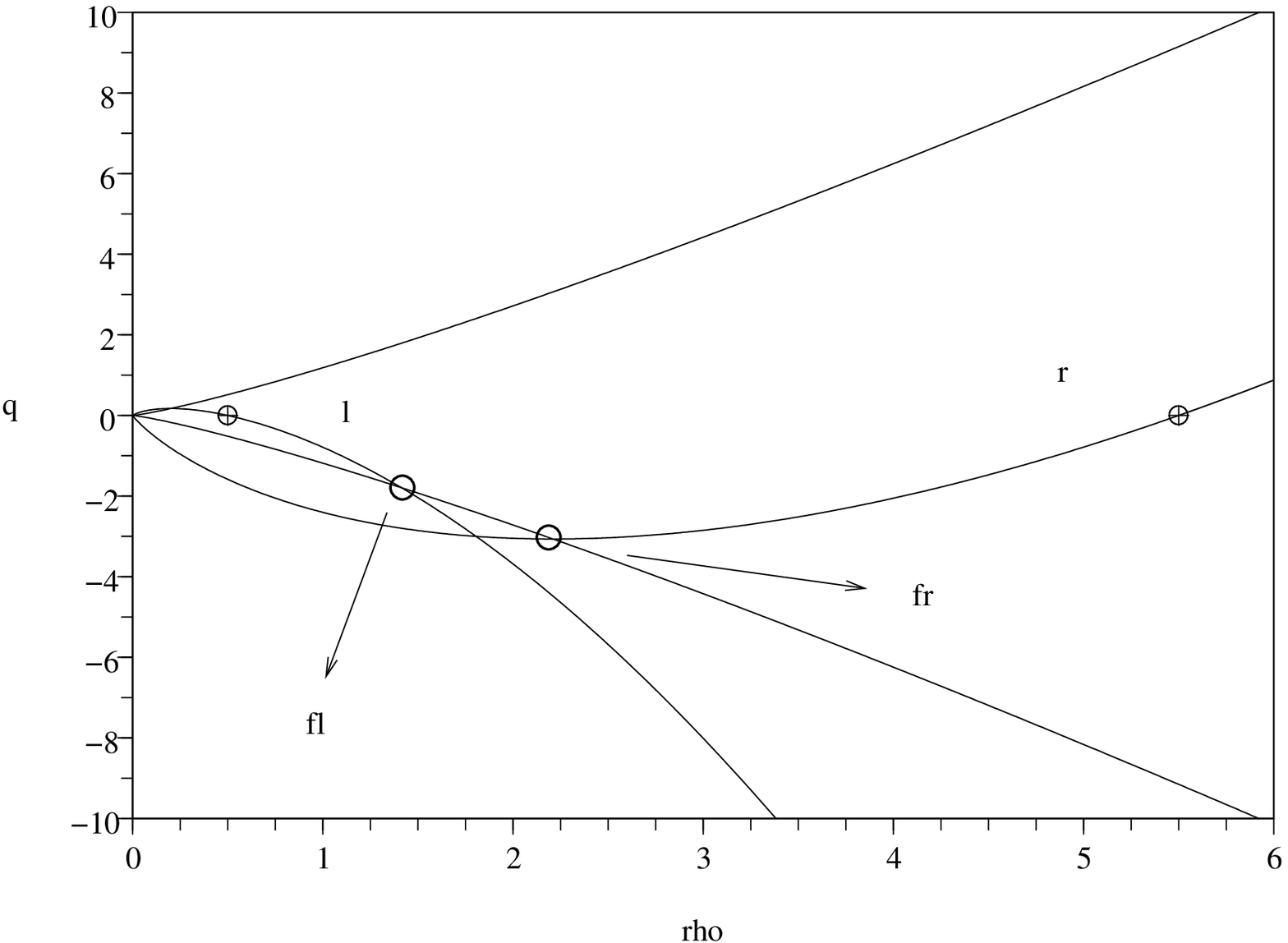}
    \end{psfrags}
  \end{minipage}
  \caption{\small{The situations of Corollaries~\ref{cor:NoL},
      \ref{cor:Nop} and~\ref{cor:NoP}.}}
  \label{fig:prop31}
\end{figure}

\begin{corollary}
  \label{cor:NoL}
  Under the assumptions above, in each of the two cases
  \begin{displaymath}
    \begin{array}{rcl@{\quad \mbox{ and } \quad}rcl}
      a_ l & < & a_r
      &
      \fd{l} (\bar u_l) & < & \fd{r} (\bar u_r)
      \\
      a_l & > & a_r
      &
      \fu{l} (\bar u_l) & > & \fu{r} (\bar u_r)
    \end{array}
  \end{displaymath}
  see Figure~\ref{fig:prop31}, an~\textbf{(L)}-solution
  to~(\ref{eq:RP}), attaining values in $A_0$, does not exist.
\end{corollary}

\par From the physical point of view, when $a_l = a_r$
condition~\textbf{(L)} is the most reasonable one, since it states the
conservation of mass and linear momentum. In the case of an elbow, an
analog of this condition can be justified through the conservation of
the linear momentum along a direction dependent on the geometry of the
elbow, see~\cite[Propositions~3.2 and~4.2]{ColomboGaravello2}. For a
study of the dynamic of a fluid in a kink, see
also~\cite{HoldenRisebroKink}. However, in the present case of a
junction between collinear pipes with different sections,
condition~\textbf{(L)} is hardly acceptable,
see~\cite[Figure~7]{ColomboGaravello2}.

%
%
%
%
\subsection{(p)-solutions}
\label{sec:p-sol}

Under condition~\textbf{(p)}, (\ref{eq:Sigma}) reads
\begin{displaymath}
  a_l \, q_l  = a_r \, q_r
  \qquad\mbox{ and }\qquad
  p(\rho_l) = p(\rho_r)
\end{displaymath}
and was considered in~\cite[Paragraph~5.2]{BandaHertyKlar1}, see
also~\cite{BandaHertyKlar2}. It was there introduced neglecting the
pressure drop at the junction, on the basis of engineering literature
on the subject, see for instance~\cite[Section
ì6.3.1]{WinterbonePearson}. In this condition, the role of the fluid
speed is limited to ensure the conservation of mass.

\begin{theorem}
  \label{th:p-sol}
  Let~\textbf{(EoS)} hold and fix positive $a_l, a_r$.  Consider the
  Riemann Problem~(\ref{eq:RP}) with $u_l, u_r \in A_0$. Define
  \begin{displaymath}
    l' = \max \left\{\fu{l}(\bar u_l),
      \fd{r}(\bar u_r) \right\} 
    \quad \textrm{ and } \quad
    l'' = \min \left\{\fd{l}(\bar u_l),
      \fu{r}(\bar u_r) \right\}. 
  \end{displaymath}
  A necessary and sufficient condition for the existence and
  uniqueness of a~\textbf{(p)}-solution to~(\ref{eq:RP}) attaining
  values in $A_0$ is
  \begin{equation}
    \label{eq:p-sol_CNS}
    \left\{
      \begin{array}{rcl}
        a_l \, L_1(l'; \bar u_l) 
        & \geq &
        a_r \, L_2^-(l'; \bar u_r),
        \\
        a_l \, L_1(l''; \bar u_l) 
        & \leq &
        a_r \, L_2^-(l''; \bar u_r). 
      \end{array}
    \right.
  \end{equation}
\end{theorem}

An example of non existence of~\textbf{(p)}-solutions to~(\ref{eq:RP})
is provided by the following corollary.

\begin{corollary}
  \label{cor:Nop}
  Under the assumptions above, in each of the two cases
  \begin{displaymath}
    \fu{l} (\bar u_l) > \fu{r} (\bar u_r)
    \quad \mbox{ or } \quad
    \fd{l} (\bar u_l) < \fd{r} (\bar u_r)
  \end{displaymath}
  a~\textbf{(p)}-solution to~(\ref{eq:RP}) attaining values in $A_0$
  does not exist.
\end{corollary}

An example of lack of continuous dependence for \textbf{(p)}-solutions
is in~\cite{ColomboGaravello2}.

>From the physical point of view, condition~\textbf{(p)} is fully
justified in the static situation $\bar q_l = \bar q_r = 0$. It is
sometimes extended to non static cases in the engineering literature,
possibly corrected through suitable \emph{pressure loss} coefficients,
see~\cite[Section~6.3.2]{WinterbonePearson}.

%
%
%
%
\subsection{(P)-solutions}
\label{sec:P-sol}

Under condition~\textbf{(P)}, (\ref{eq:Sigma}) reads
\begin{displaymath}
  a_l \, q_l  = a_r \, q_r
  \qquad\mbox{ and }\qquad
  P(\rho_l, q_l) = P(\rho_r, q_r)
\end{displaymath}
and was considered in~\cite{ColomboGaravello1, ColomboGaravello3}.

\begin{theorem}
  \label{th:P-sol}
  Let~\textbf{(EoS)} hold and fix positive $a_l, a_r$.  Consider the
  Riemann Problem~(\ref{eq:RP}) with $\bar u_l, \bar u_r \in
  A_0$. Define
  \begin{eqnarray*}
    l'
    & = &
    \inf \left\{ \rho \in [\fu{l}(\bar u_l),
      \fd{l}(\bar u_l)] \colon
      a_l \, L_1(\rho; \bar u_l) \le
      a_r \, L_2^-(\fu{r}(\bar u_r); \bar u_r)
    \right\}
    \\
    l''
    & = &
    \sup \left\{ \rho \in [\fu{l}(\bar u_l),
      \fd{l}(\bar u_l)] \colon
      a_l \, L_1(\rho; \bar u_l) \ge
      a_r \, L_2^-(\fd{r}(\bar u_r); \bar u_r)
    \right\}
  \end{eqnarray*}
  A necessary and sufficient condition for the existence and
  uniqueness of a~\textbf{(P)}-solution to~(\ref{eq:RP}), attaining
  values in $A_0$, is
  \begin{equation}
    \label{eq:P-sol_CNS}
    \left\{
      \begin{array}{rcl}
        P \left(l', L_1(l'; \bar u_l)\right) - 
        P\left (g \left(\frac{a_l}{a_r} L_1(l'; \bar u_l) \right),
          \frac{a_l}{a_r} L_1(l'; \bar u_l)\right) 
        & \leq & 0,
        \\
        P \left(l'', L_1(l''; \bar u_l)\right) -  
        P
        \left( g \left(\frac{a_l}{a_r} L_1(l''; \bar u_l)\right),
          \frac{a_l}{a_r} L_1(l''; \bar u_l)
        \right) 
        & \geq & 0,
      \end{array}
    \right.
  \end{equation}
  where $g$ is the inverse in $A_0$ of the function $\rho \to
  L_2^-(\rho ; \bar u_r)$.
\end{theorem}

The proof is similar to that of Theorem~\ref{thm:Lsol}; hence we omit
it.  An example of non existence of~\textbf{(P)}-solutions
to~(\ref{eq:RP}) is provided by the following corollary.

\begin{corollary}
  \label{cor:NoP}
  Under the assumptions above, in each of the two cases
  \begin{displaymath}
    \fu{l} (\bar u_l) > \fu{r} (\bar u_r)
    \quad \mbox{ or } \quad
    \fd{l} (\bar u_l) < \fd{r} (\bar u_r),
  \end{displaymath}
  see Figure~\ref{fig:prop31}, a~\textbf{(P)}-solution
  to~(\ref{eq:RP}) attaining values in $A_0$ does not exist.
\end{corollary}

>From the physical point of view, condition~\textbf{(P)} implies the
conservation of linear momentum along directions orthogonal to the
pipes, as shown in~\cite[Lemma~2.2]{ColomboGaravello3}.

%
%
%
%
\subsection{(S)-solutions}
\label{sec:s-sol}

Under condition~\textbf{(S)}, (\ref{eq:Sigma}) reads
\begin{displaymath}
  a_l \, q_l  = a_r \, q_r
  \qquad\mbox{ and }\qquad
  a_r P(\rho_r, q_r) = a_l P(\rho_l, q_l) + \int_{a_l}^{a_r}
  p\left( R(\alpha; \rho_l, q_l) \right) d\alpha,
\end{displaymath}
where $\left( R(\alpha; \rho_l, q_l), Q(\alpha; \rho_l, q_l) \right)$
is the solution to the Cauchy Problem
\begin{equation}
  \label{eq:sysR}
  \left\{
    \begin{array}{l}
      \frac{d}{da} Q(a) = - Q(a) / a,
      \vspace{.2cm}
      \\
      \frac{d}{da} \left[ a P \left(R(a), Q(a) \right) \right] 
      = p \left(R(a) \right),
      \vspace{.2cm}
      \\
      Q(a_l) = q_l,
      \vspace{.2cm}
      \\
      R(a_l) = \rho_l.
    \end{array}
  \right.
\end{equation}
This kind of solution was considered in~\cite{ColomboMarcellini} as
limit of solutions to the $p$-system with smoothly varying
section. This argument is based on the next lemma.

\begin{lemma}
  \label{lem:Sigma}
  Fix $\check x < \hat x \in \reali$, $a_l>0$ and $a_r>0$ with $a_l
  \neq a_r$. Let $a \in \C{0,1} (\reali;\pint{\reali}^+ )$ satisfy
  \begin{equation}
    \label{eq:CondA}
    \left\{
      \begin{array}{l@{\quad\mbox{ if } \quad}rcl}
        a(x) = a_l & x & < & \check x\\
        a \mbox{ strictly monotone} & x & \in & [\check x, \hat x]\\
        a(x) = a_r & x & > & \hat x
      \end{array}
    \right.
  \end{equation}
  Call $\tilde R^a (x; \rho_l, q_l )$ the $\rho$ component of the
  solution to the Cauchy Problem
  \begin{displaymath}
    \left\{
      \begin{array}{l@{\qquad}rcl}
        \partial_x (a(x) \,  q) = 0
        &
        \rho(\check x) & = & \rho_l
        \\
        \partial_x \!
        \left( 
          \!
          a(x) \,\left( q^{2} / \rho + p(\rho) \right)
          \!
        \right) 
        =
        p(\rho)\, \partial_x a
        &
        q(\check x) & = & q_l \,.
      \end{array}
    \right.
  \end{displaymath}
  Then, the function
  \begin{equation}
    \label{eq:SigmaSmooth}
    (a_l,a_r;\rho_l,q_l; \check x, \hat x; a) \to
    \int_{\check x}^{\hat x} p \left( R^a(x;\rho_l,q_l) \right) a'(x) \, dx \,.
  \end{equation}
  is well defined. Moreover, if $a^1, a^2$ are monotone functions
  satisfying~(\ref{eq:CondA}), then the corresponding
  functions~(\ref{eq:SigmaSmooth}) coincide. Hence, the following map
  is well defined:
  \begin{equation}
    \label{eq:SigmaAlfa}
    (a_l,a_r;\rho_l,\rho_r) \to
    \int_{a_l}^{a_r} p \left( R(\alpha;\rho_l,q_l)\right) \, d\alpha \,.
  \end{equation}
\end{lemma}

For a proof of this result,
see~\cite[Proposition~2.7]{ColomboMarcellini}.  Note that the Cauchy
Problem~(\ref{eq:sysR}) can be rewritten in the form
\begin{equation}
  \label{eq:sysRQ_explicit}
  \left\{
    \begin{array}{l}
      Q(a) = \frac{a_l}{a}\, q_l,\vspace{.2cm}
      \\
      \frac{d}{da} R(a) = \frac{R(a)}{a}\,
      \frac{a_l^2 q_l^2}{a^2 p'\left(R(a) \right) 
        R^2(a) - a_l^2 q_l^2} \vspace{.2cm}
      \\
      R(a_l) = \rho_l.
    \end{array}
  \right.
\end{equation}

\begin{figure}[htpb]
  \centering
  \begin{psfrags}
    \psfrag{q}{$q$} \psfrag{rho}{$\rho$}
    \includegraphics[width=6cm]{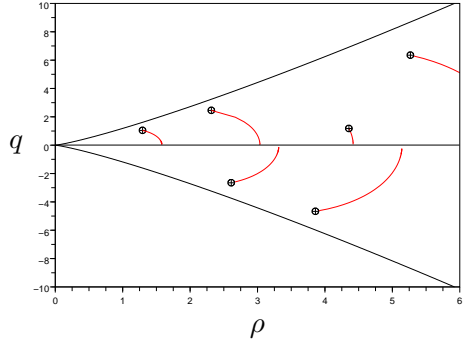}
  \end{psfrags}
  \caption{\small{Solutions to~(\ref{eq:sysRQ_explicit}) starting from
      various initial conditions.}}
  \label{fig:S_curves}
\end{figure}

\begin{lemma}
  \label{le:5.1}
  Consider the system~(\ref{eq:sysRQ_explicit}) and a point $(\rho_l,
  q_l) \in A_0$.
  \begin{enumerate}
  \item If $q_l > 0$, then $Q(a)$ is strictly decreasing.  If $q_l =
    0$, then $Q(a)$ is constantly equal to $0$.  If $q_l < 0$, then
    $Q(a)$ is strictly increasing.

  \item If $q_l = 0$, then $R(a)$ is constantly equal to $\rho_l$.  If
    $q_l \ne 0$, then $R(a)$ exists and is strictly increasing for
    every $a > a_l$. Moreover it is bounded for $a > a_l$; see
    Figure~\ref{fig:S_curves}.
  \end{enumerate}
\end{lemma}

\begin{theorem}
  \label{th:S-sol}
  Let~\textbf{(EoS)} hold and fix positive $a_l < a_r$.  Consider the
  Riemann Problem~(\ref{eq:RP}) with $u_l, u_r \in A_0$.  Define
  \begin{eqnarray*}
    l'
    & = &
    R \left( a_r; \fu{l} (\bar u_l),
      L_1 \left( \fu{l} (\bar u_l); \bar u_l
      \right) \right),
    \\
    l''
    & = &
    R \left( a_r; \fd{l} (\bar u_l),
      L_1 \left( \fd{l} (\bar u_l); \bar u_l
      \right) \right),
  \end{eqnarray*}
  where $R$ is the function defined in~(\ref{eq:sysRQ_explicit}).  A
  necessary and sufficient condition for existence and uniqueness of
  an~\textbf{(S)}-solution to~(\ref{eq:RP}), attaining values in
  $A_0$, is
  \begin{equation}
    \label{eq:S-sol_CNS}
    \left\{
      \begin{array}{rcl}
        a_l \, L_1 \left( \fu{l} (\bar u_l); \bar u_l \right) 
        & \geq &
        a_r \, L_2^- \left( l'; \bar u_r \right),
        \\
        a_l \, L_1 \left( \fd{l} (\bar u_l); \bar u_l \right) 
        & \leq &
        a_r \, L_2^- \left( l''; \bar u_r \right).
      \end{array}
    \right.
  \end{equation}
\end{theorem}

\begin{corollary}
  \label{cor:S_noexistence}
  Consider the Riemann Problem~(\ref{eq:RP}) with $a_l < a_r$ and such
  that $\fu{l} (\bar u_l) > \fu{r} (\bar u_r)$.  Then
  an~\textbf{(S)}-solution to~(\ref{eq:RP}), attaining values in
  $A_0$, does not exist.
\end{corollary}

>From the physical point of view, condition~\textbf{(S)} is justified
as the limit of smooth changes in the pipes' section,
see~\cite[Theorem~2]{GuerraMarcelliniSchleper}. Viceversa, consider
$n$ consecutive junctions sited at, say, $x^n_i = i/n$, with $i=0, 1 ,
\ldots, n$, separating pipes with section, say, $a^n_i = \bar a_l +
(i/n) (\bar a_r - \bar a_l)$. Imposing condition~\textbf{(S)} on each
junction is equivalent, in the limit $n \to +\infty$, to the usual
model for pipes with a smoothly varying section,
see~\cite[Section~8.1]{Whitham}.

\Section{Numerical Examples}
\label{sec:Numerical}

The paragraph is devoted to show results of numerical integrations of
solutions to~(\ref{eq:System})--(\ref{eq:Phi}) with the different
choices of the function $\Phi$ considered above. The integrations
below are \emph{``exact''}, in the sense that they amount to the
solutions of Riemann problems at the junction, which are obtained
through the (approximate) computation of the intersection between Lax
curves~(\ref{eq:LaxCurves}). Neither time steps, nor meshes, are
involved.

Throughout, we fix $\gamma = 1.4$ and the $\gamma$-law $p(\rho) =
\rho^\gamma$.

\begin{figure}[htpb]
  \centering
  \begin{psfrags}
    \psfrag{rho}{{\tiny $\rho$}} \psfrag{q}{{\tiny $q$}}
    \psfrag{x}{{\tiny $x$}}%
    \psfrag{L-solution - q}{{\tiny\textbf{(L)}-solution - $q$}}%
    \psfrag{L-solution - rho}{{\tiny\textbf{(L)}-solution - $\rho$}}%
    \psfrag{p-solution - q}{{\tiny \textbf{(p)}-solution - $q$}}%
    \psfrag{p-solution - rho}{{\tiny \textbf{(p)}-solution - $\rho$}}%
    \psfrag{P-solution - q}{{\tiny \textbf{(P)}-solution - $q$}}%
    \psfrag{P-solution - rho}{{\tiny \textbf{(P)}-solution - $\rho$}}%
    \psfrag{S-solution - q}{{\tiny \textbf{(S)}-solution - $q$}}%
    \psfrag{S-solution - rho}{{\tiny \textbf{(S)}-solution - $\rho$}}%
    \includegraphics[width=6cm]{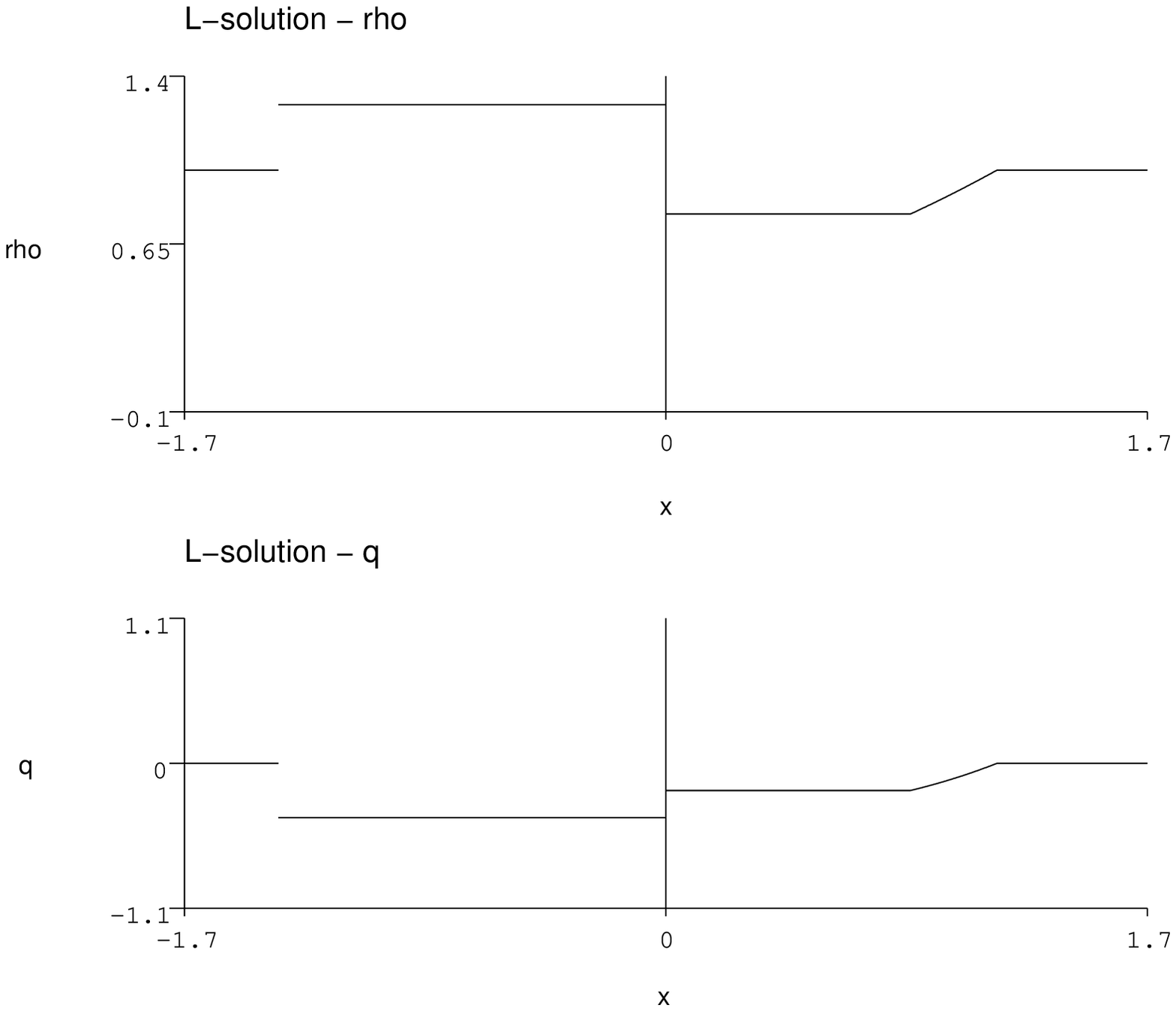}
    \includegraphics[width=6cm]{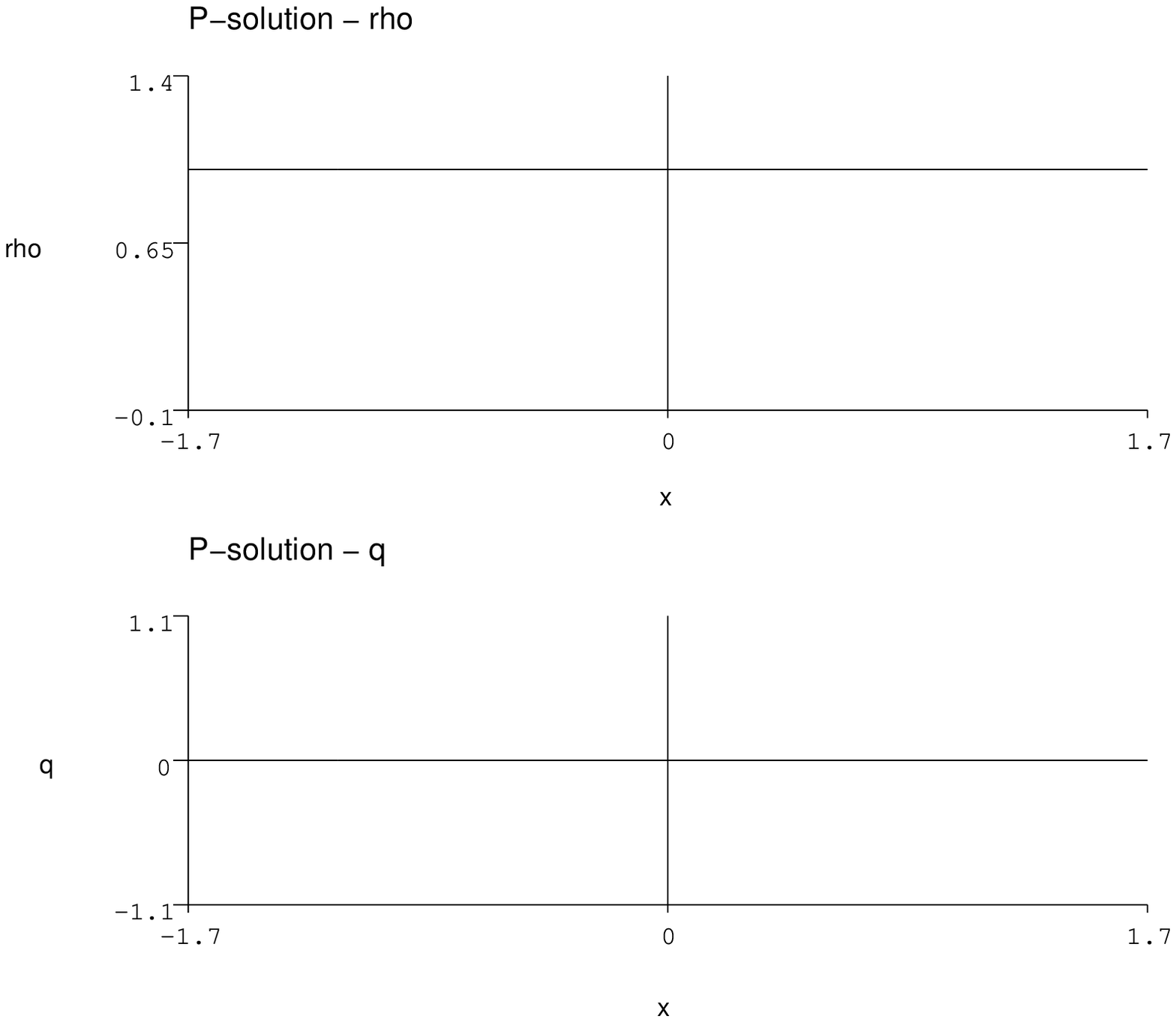}\\
    \includegraphics[width=6cm]{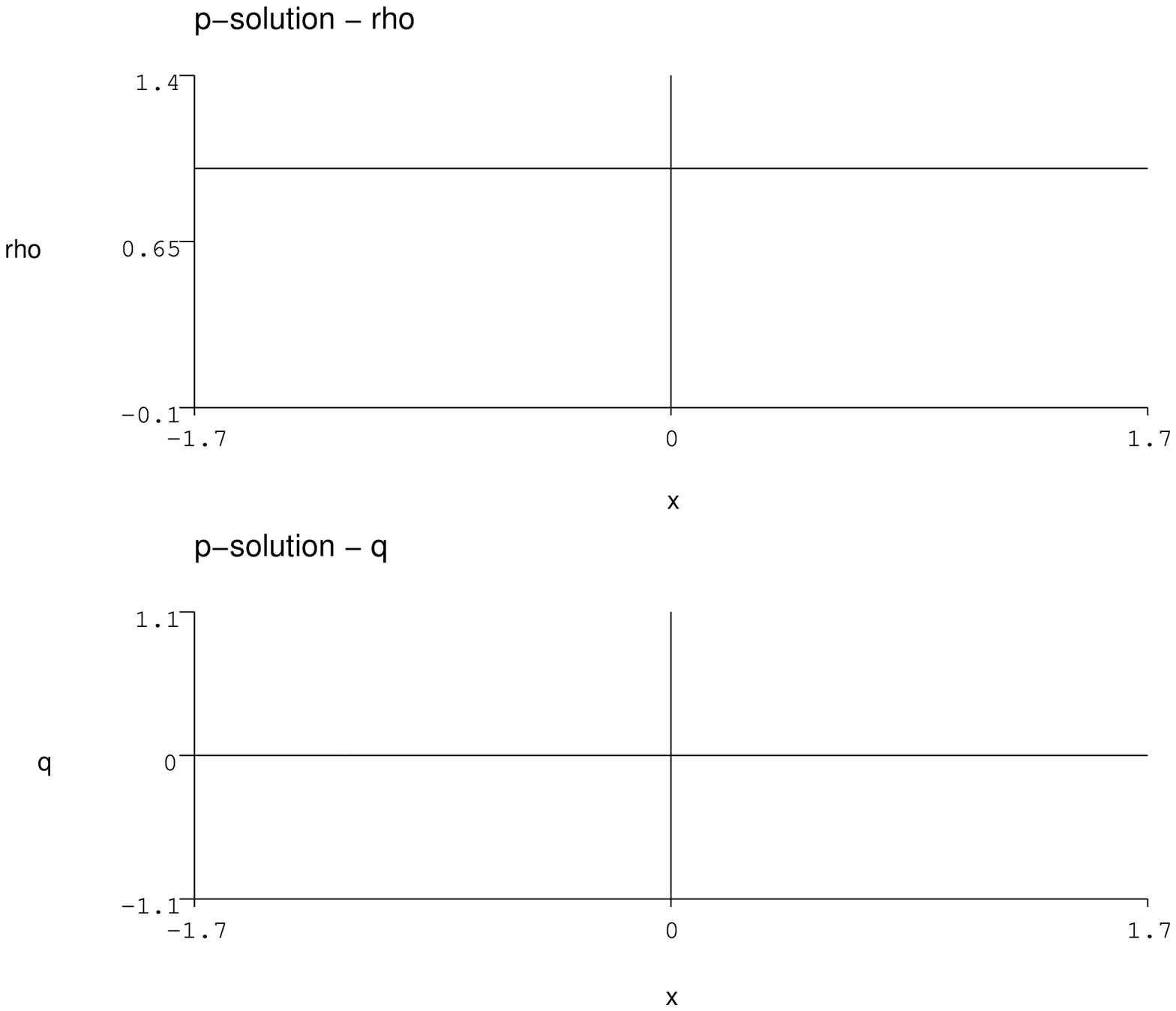}
    \includegraphics[width=6cm]{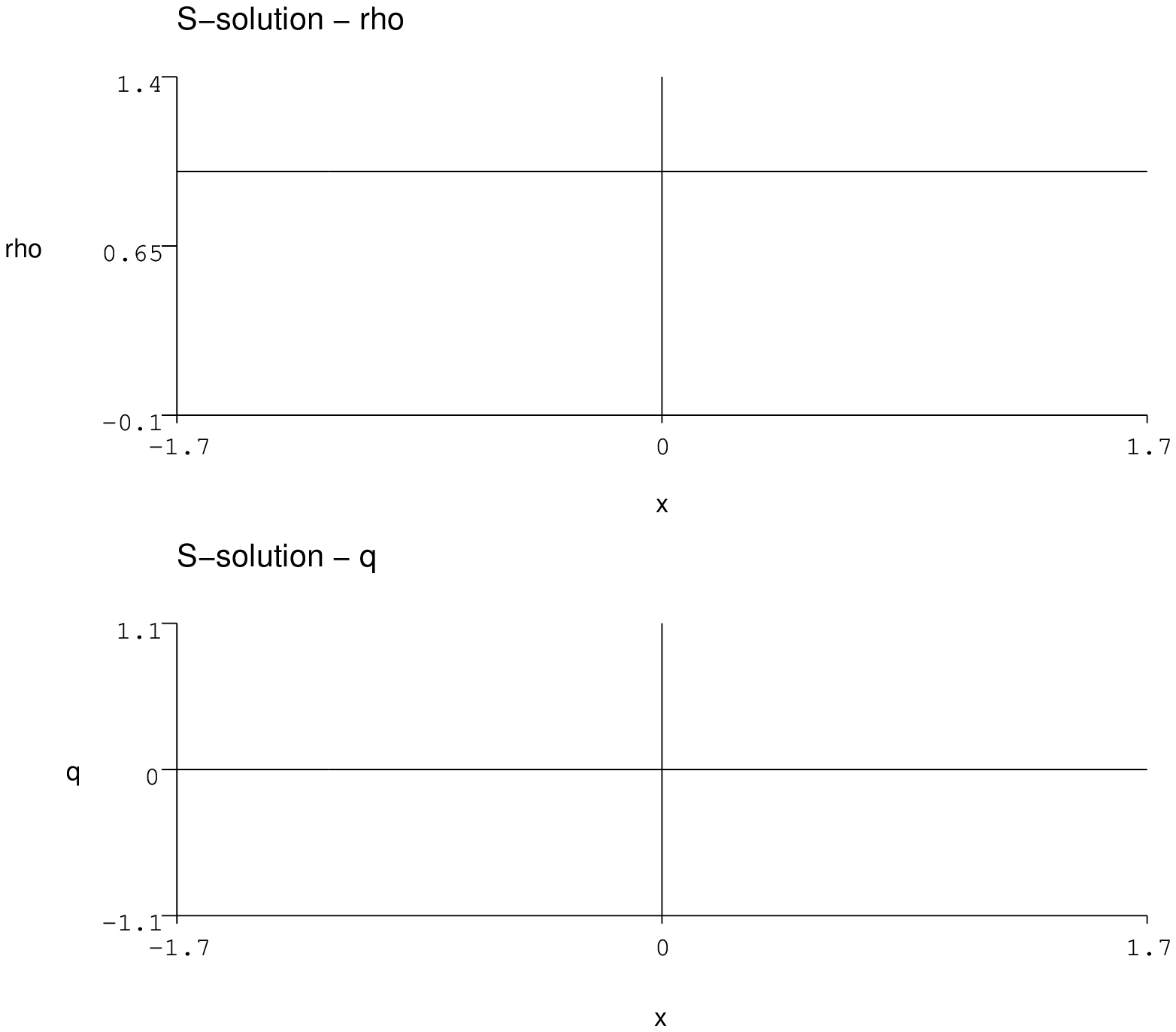}
  \end{psfrags}
  \caption{\small{Here, $a_l=1$, $a_r=2$, $\bar\rho_l = \bar\rho_r=1$,
      $\bar q_l =\bar q_r = 0$. Note that only
      the~\textbf{(L)}-solution is \emph{not} stationary.}}
  \label{fig:caso1}
\end{figure}
First, consider an initial datum with the same density and zero speed
in the two pipes: $\bar q_l = \bar q_r = 0$ and $\bar \rho_l = \bar
\rho_r$. Conditions~\textbf{(P)}, \textbf{(p)} and~\textbf{(S)} yield
the stationary solution: no transfer of fluid between the pipes and no
pressure difference at the junction, see Figure~\ref{fig:caso1}.  The
\textbf{(L)}-solution significantly differs: it is not stationary,
prescribes a transfer of fluid from the right (larger) tube to the
left (smaller) tube and yields a pressure difference at the junction,
see Figure~\ref{fig:caso1}. Indeed, the~\textbf{(L)}-solution consists
in a $2$-rarefaction moving to the right and a $1$-shock to the left,
with fluid flowing from the right to the left.

A situation frequently considered in the engineering literature, see
for instance~\cite[Section~6.1.3]{WinterbonePearson}
and~\cite{Benson}, is that of a shock wave hitting a junction.
\begin{figure}[htpb]
  \centering
  \begin{psfrags}
    \psfrag{rho}{{\tiny $\rho$}} \psfrag{q}{{\tiny $q$}}
    \psfrag{x}{{\tiny $x$}}%
    \psfrag{L-solution - q}{{\tiny\textbf{(L)}-solution - $q$}}%
    \psfrag{L-solution - rho}{{\tiny\textbf{(L)}-solution - $\rho$}}%
    \psfrag{p-solution - q}{{\tiny \textbf{(p)}-solution - $q$}}%
    \psfrag{p-solution - rho}{{\tiny \textbf{(p)}-solution - $\rho$}}%
    \psfrag{P-solution - q}{{\tiny \textbf{(P)}-solution - $q$}}%
    \psfrag{P-solution - rho}{{\tiny \textbf{(P)}-solution - $\rho$}}%
    \psfrag{S-solution - q}{{\tiny \textbf{(S)}-solution - $q$}}%
    \psfrag{S-solution - rho}{{\tiny \textbf{(S)}-solution - $\rho$}}%
    \includegraphics[width=6cm]{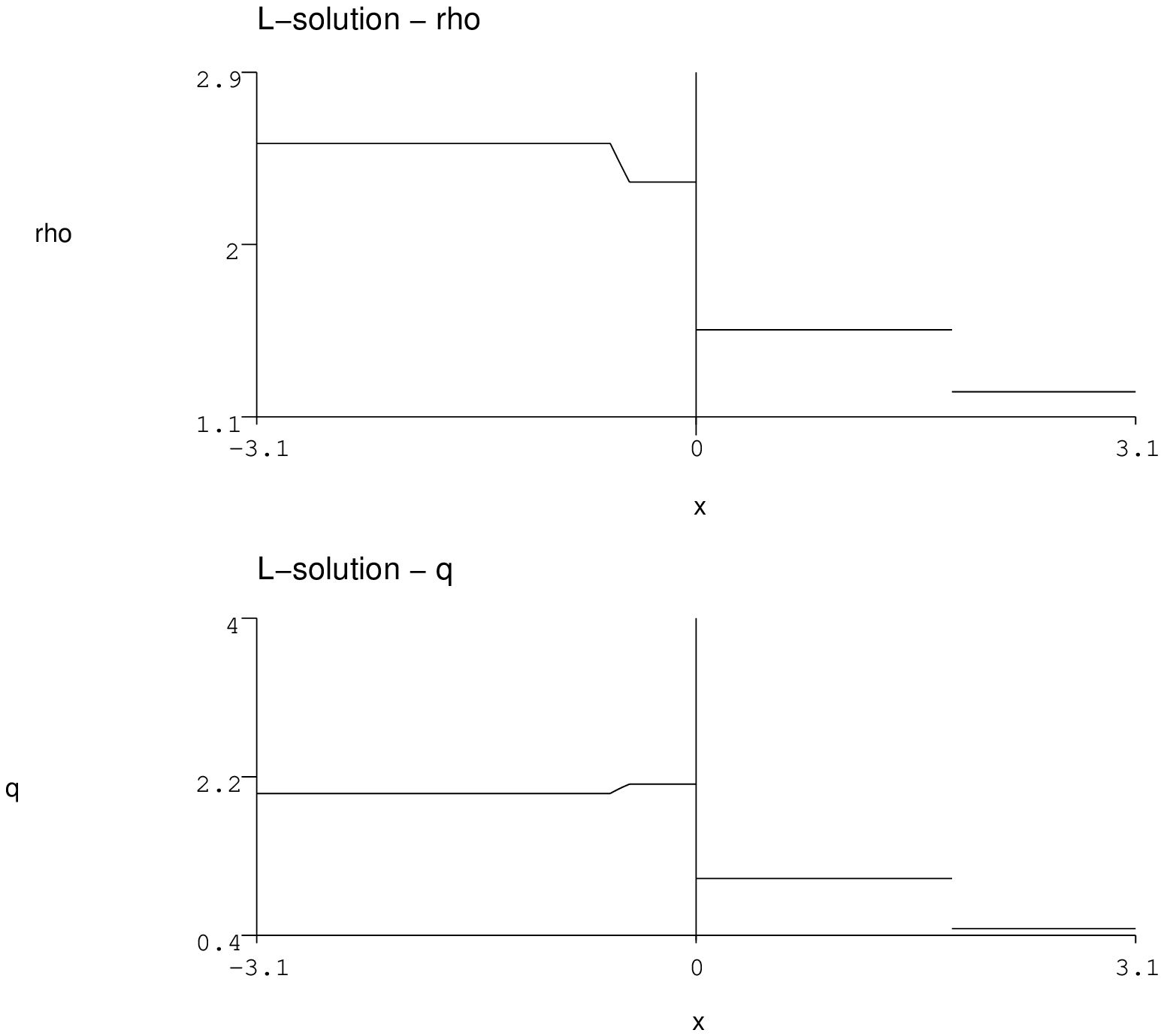}
    \includegraphics[width=6cm]{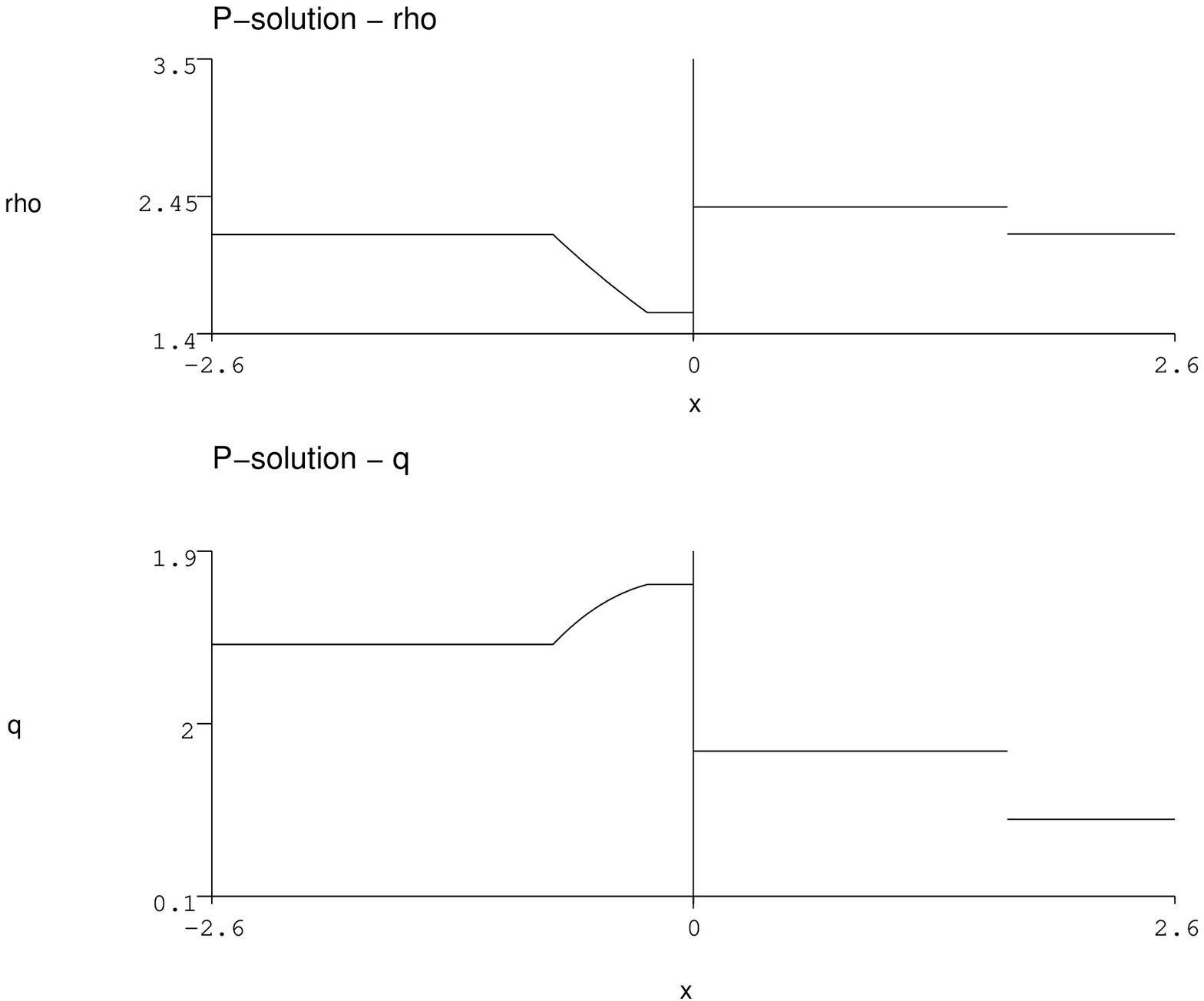}\\
    \includegraphics[width=6cm]{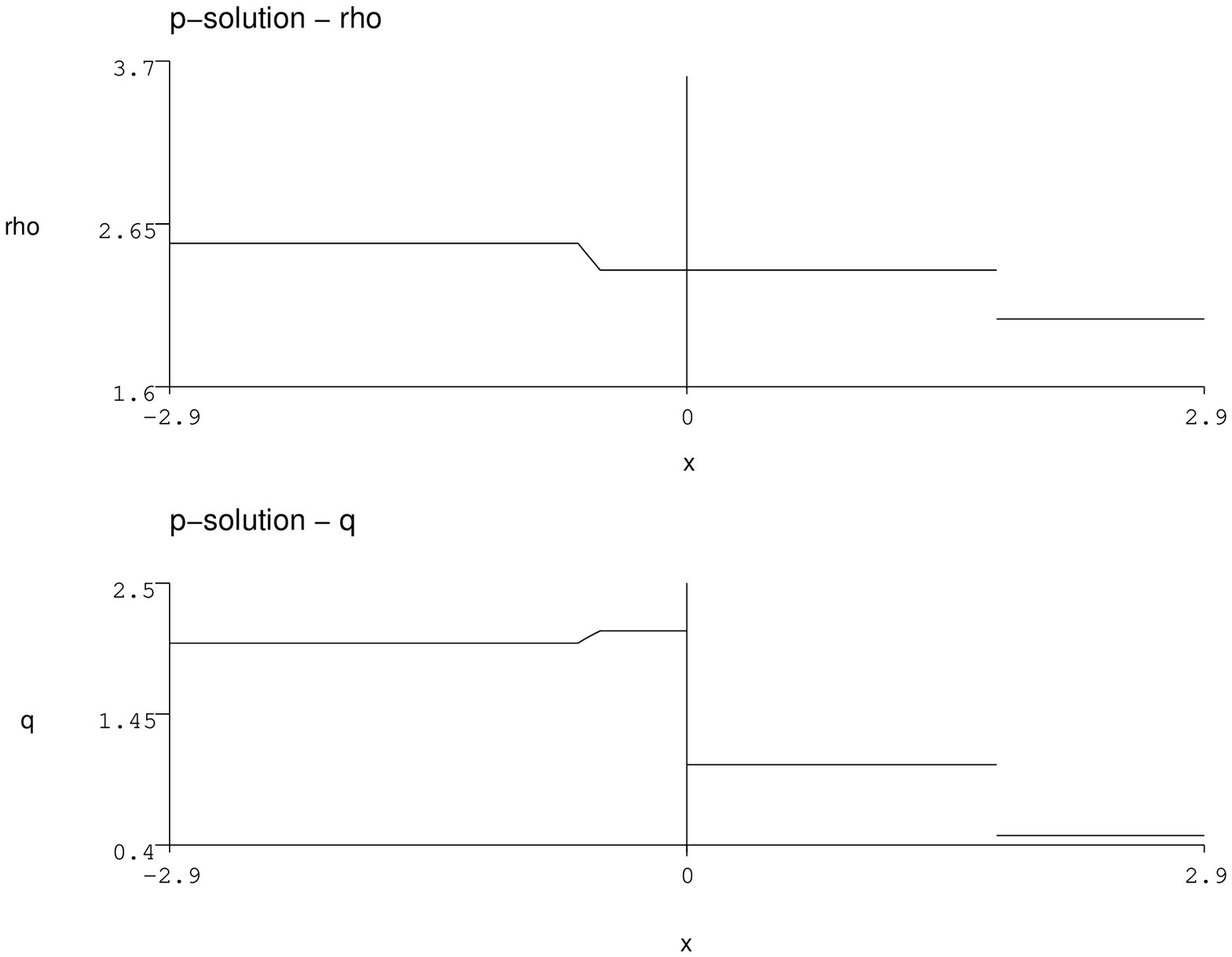}
    \includegraphics[width=6cm]{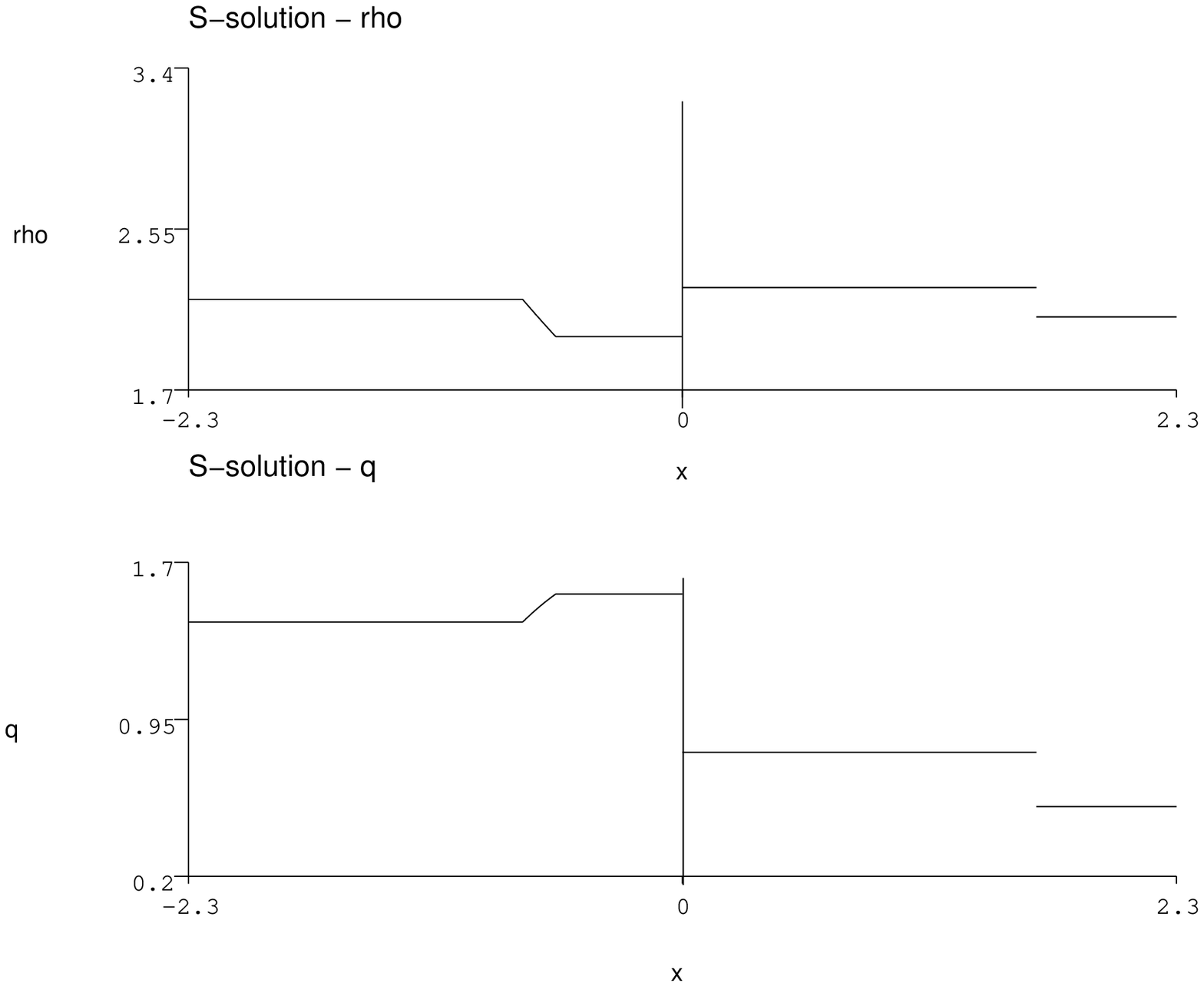}
  \end{psfrags}
  \caption{\small{Here, $a_l=1$, $a_r=2$, and the left state is
      $\bar\rho_l = 2.2$, $\bar q_l = 1.387882$. For the other initial
      data see~(\ref{eq:table}).}}
  \label{fig:sa}
\end{figure}
We consider the case of stationary flow of fluid moving from the
smaller pipe, say, on the left, to the larger one on the
right. Differently from the previous case, we now want to perturb a
stationary but non-static situation. Therefore, we are bound to choose
different stationary configurations for the different solutions. Let
$a_l = 1$, $a_r = 2$ and fix the unperturbed state $\rho_l=2$, $q_l =
1$ on the left. Note that $q_l > 0$, so that fluid flows
rightwards. Correspondingly, we find the state $(\bar\rho_r, \bar
q_r)$ to be assigned to the right pipe, so that we have stationary
solutions in the different cases:
\begin{equation}
  \label{eq:table}
  \begin{array}{cll}
    \mbox{Solution} & \bar\rho_r & \bar q_r
    \\
    \textbf{(L)} & 1.2520452 & 0.5
    \\
    \textbf{(P)} & 2.2051669 & 0.5
  \end{array}
  \qquad
  \begin{array}{cll}
    \mbox{Solution} & \bar\rho_r & \bar q_r
    \\
    \textbf{(p)} & 2.0000000 & 0.5
    \\
    \textbf{(S)} & 2.1064869 & 0.5
  \end{array}
\end{equation}
Obviously, $\bar q_r$ is uniquely determined by $q_l$ and the
sections $a_l$, $a_r$ through mass conservation. Then, a $2$-shock
moving in the left tube towards the junction has left state
$\bar\rho_l = 2.2$, $\bar q_l = 1.387882$ and right state $(\rho_l,
q_l)$. These latter values, as well as the forthcoming solutions, are
determined using~(\ref{eq:LaxCurves}) and the definitions of
solutions.

The qualitative behavior of the solution to the Riemann
Problem~(\ref{eq:RP}) is the same in all cases: the shock interacts
with the junction leading to the formation of two waves. A refracted
shock proceeding in the right tube and a reflected rarefaction moving
leftwards, see~Figure~\ref{fig:sa}.

Consider now a stationary configuration perturbed by a shock coming
from the right (larger) pipe.
\begin{figure}[htpb]
  \centering
  \begin{psfrags}
    \psfrag{rho}{{\tiny $\rho$}} \psfrag{q}{{\tiny $q$}}
    \psfrag{x}{{\tiny $x$}}%
    \psfrag{L-solution - q}{{\tiny\textbf{(L)}-solution - $q$}}%
    \psfrag{L-solution - rho}{{\tiny\textbf{(L)}-solution - $\rho$}}%
    \psfrag{p-solution - q}{{\tiny \textbf{(p)}-solution - $q$}}%
    \psfrag{p-solution - rho}{{\tiny \textbf{(p)}-solution - $\rho$}}%
    \psfrag{P-solution - q}{{\tiny \textbf{(P)}-solution - $q$}}%
    \psfrag{P-solution - rho}{{\tiny \textbf{(P)}-solution - $\rho$}}%
    \psfrag{S-solution - q}{{\tiny \textbf{(S)}-solution - $q$}}%
    \psfrag{S-solution - rho}{{\tiny \textbf{(S)}-solution - $\rho$}}%
    \includegraphics[width=6cm]{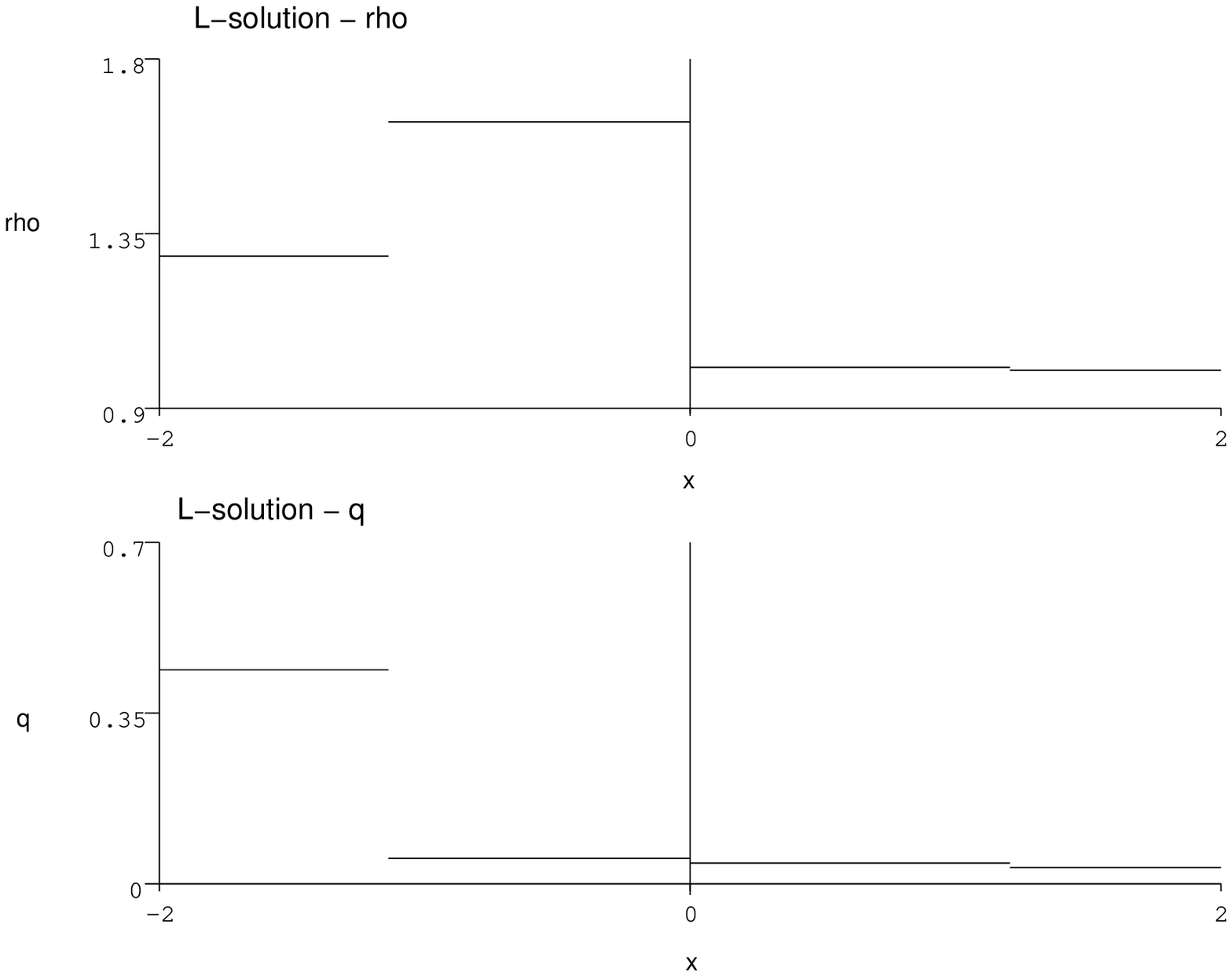}
    \includegraphics[width=6cm]{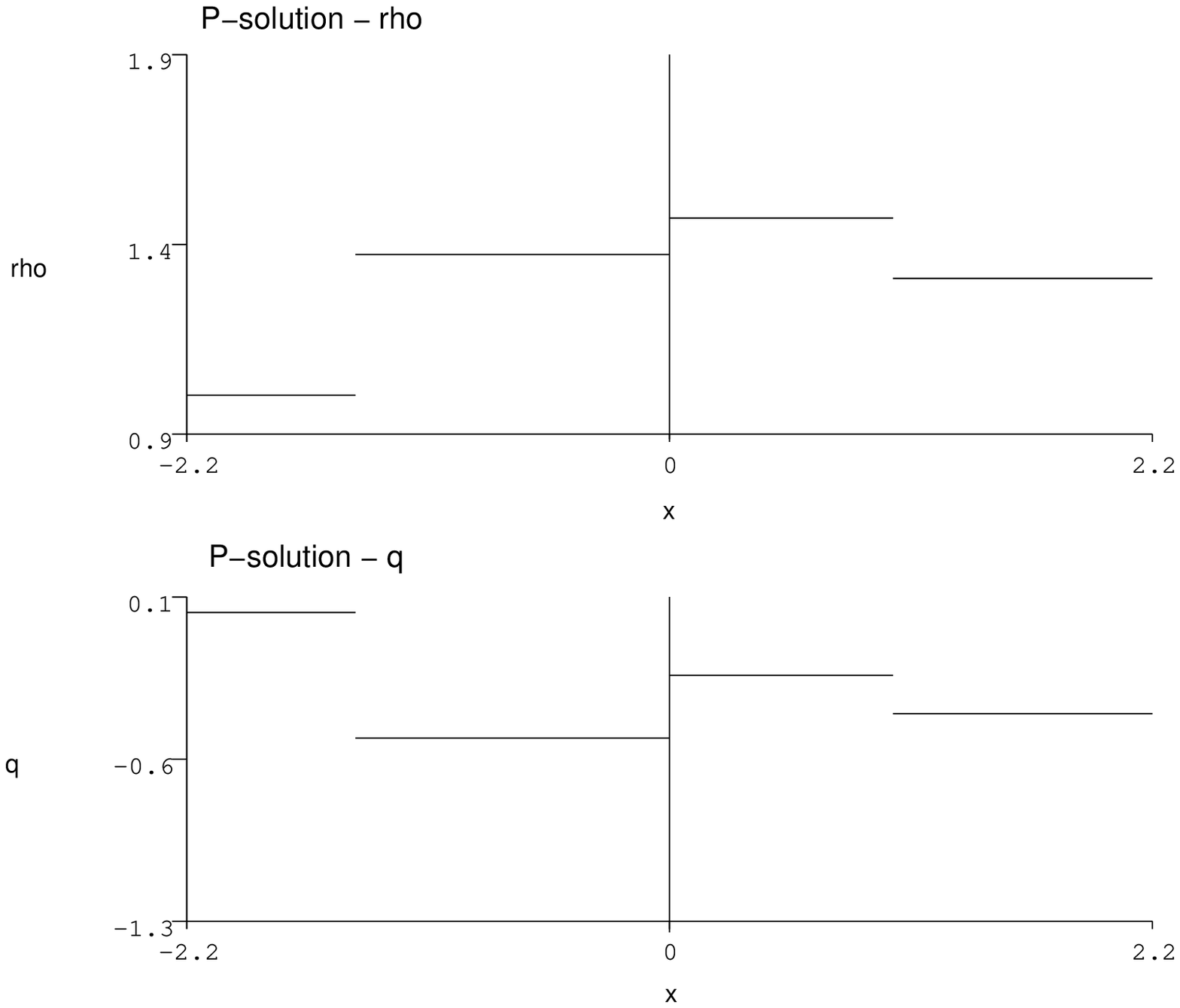}\\
    \includegraphics[width=6cm]{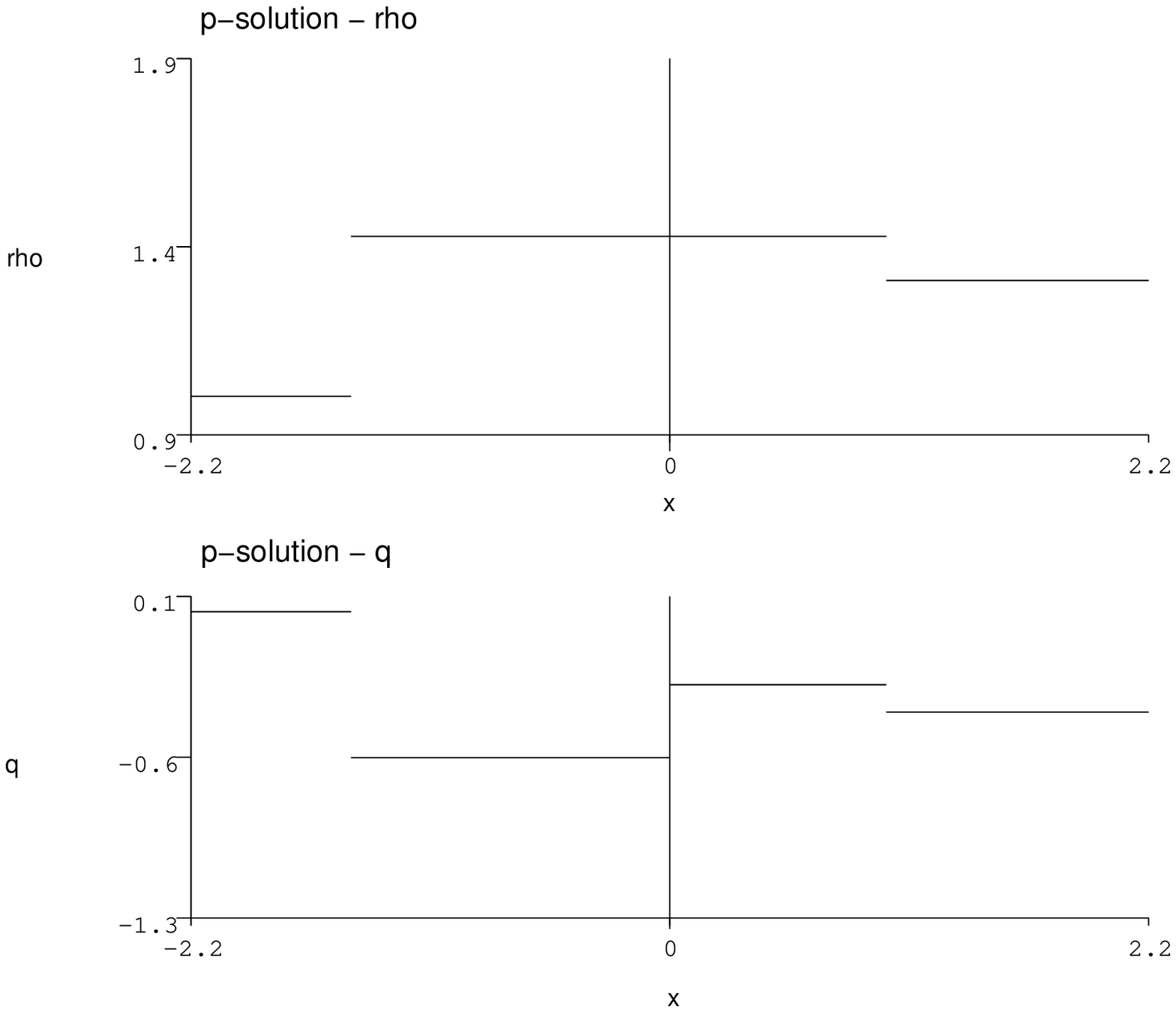}
    \includegraphics[width=6cm]{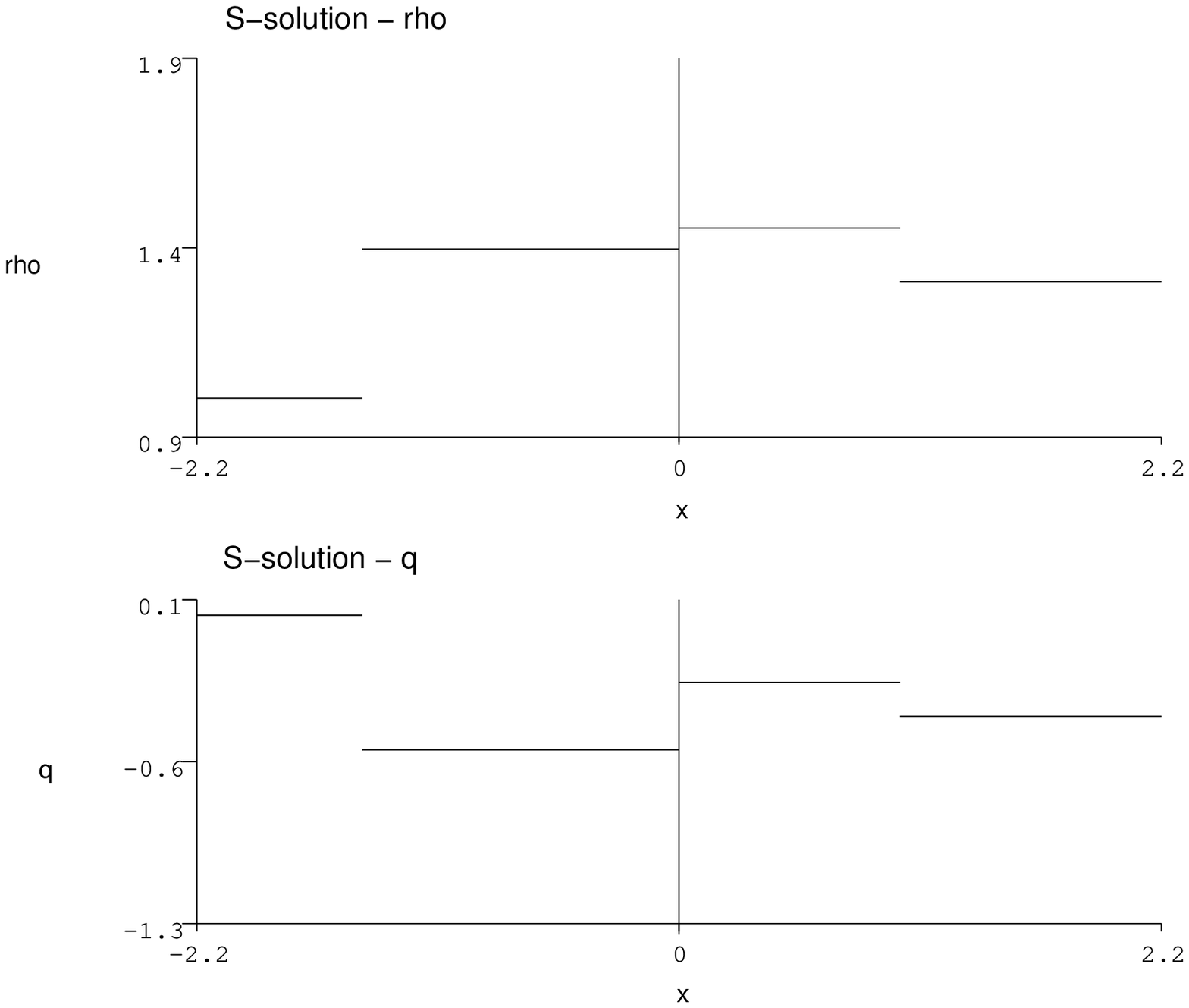}
  \end{psfrags}
  \caption{\small{Here, $a_l=1$, $a_r=2$. In case~\textbf{(L)} $\bar
      \rho_l = 2$, $\bar q_l = 1$, $\rho_r = 1.2520452$, $q_r = 0.5$.
      In cases~\textbf{(p)}, \textbf{(P)} and~\textbf{(S)}
      $\bar\rho_l = \rho_r = 1$, $\bar q_l = q_r = 0$.
      The interaction of
      a $1$--shock from the right results in a $1$ shock being
      refracted in the left tube and a $2$-shock reflected in the
      right tube.}}
  \label{fig:caso5}
\end{figure}
In cases~\textbf{(p)}, \textbf{(P)} and~\textbf{(S)} we use the same
data as in Figure~\ref{fig:caso1}, namely $a_l=1$, $a_r=2$,
$\bar\rho_l = \rho_r = 1$, $\bar q_l = q_r = 0$, so that the
fluid is at rest. In case~\textbf{(L)}, the data $a_l = 1$, $a_r = 2$,
$\bar \rho_l =2$, $\bar q_l = 1$, $\rho_r = 1.2520452$ and $q_r = 0.5$
yield a stationary solution. We perturb these stationary solutions
with a $1$-shock heading
towards the junction from the right tube, obtaining the Riemann
problem with right data
$\bar\rho_r = 1.3$ and $\bar q_r = - 0.4160525$
in cases~\textbf{(p)}, \textbf{(P)} and~\textbf{(S)},
and with right data $\bar \rho_r = 1.55$ and $\bar q_r =
0.1994594$ in case~\textbf{(L)}.
Then, the result of this interaction is in
Figure~\ref{fig:caso5}. In all cases, the interaction results in two
shocks, one reflected and one refracted.

%
%
%
%
\Section{Technical Details}
\label{sec:Tech}

Recall the expressions of the eigenvalues and eigenvectors of the
$p$-system
\begin{equation}
  \label{eigenvectors}
  \begin{array}{rcl@{\quad}rcl}
    \lambda_{1} (u)
    & = &
    \frac{q}{\rho }-\sqrt{p^{\prime } (\rho)} \,,
    &
    \lambda_{2} (u)
    & = &
    \frac{q}{\rho }+\sqrt{p^{\prime } (\rho)} \,,
    \\[5pt]
    r_{1} (u)
    & = &
    \left[ 
      \begin{array}{c}
        -1 
        \\ 
        -\lambda_{1} (u)
      \end{array}
    \right] \,,
    &
    r_{2} (u)
    & = &
    \left[ 
      \begin{array}{c}
        1 
        \\ 
        \lambda_{2} (u)
      \end{array}%
    \right]\,.
  \end{array}
\end{equation}
The Lax curves and their reversed counterparts have the expressions
\begin{equation}
  \label{eq:LaxCurves}
  \begin{array}{rcl}
    L_1(\rho;u_o) 
    & = &
    \left\{
      \begin{array}{l@{\quad\mbox{ if }}l}
        \frac{\rho}{\rho_o} q_o - 
        \rho \int_{\rho_o}^{\rho} \frac{c(r)}{r}\,dr
        & \rho \leq \rho_o
        \\
        \frac{\rho}{\rho_o} q_o -
        \sqrt{\frac{\rho}{\rho_o} (\rho-\rho_o) 
          \left(p(\rho)-p(\rho_o) \right)}
        & \rho \ge \rho_o  
      \end{array}
    \right.
    \\
    L_2(\rho;u_o)
    & = &
    \left\{
      \begin{array}{l@{\quad\mbox{ if }}l}
        \frac{\rho}{\rho_o} q_o -
        \sqrt{\frac{\rho}{\rho_o} (\rho-\rho_o) 
          \left(p(\rho)-p(\rho_o) \right)}
        & \rho \leq \rho_o
        \\
        \frac{\rho}{\rho_o} q_o +
        \rho \int_{\rho_o}^{\rho} \frac{c(r)}{r}\,dr
        & \rho \ge \rho_o  
      \end{array}
    \right.
    \\
    L_1^-(\rho;u_o)
    & = &
    \left\{
      \begin{array}{l@{\quad\mbox{ if }}l}
        \frac{\rho}{\rho_o}q_o +
        \sqrt{\frac{\rho}{\rho_o} (\rho-\rho_o) 
          \left(p(\rho)-p(\rho_o) \right)}
        & \rho \leq \rho_o
        \\
        \frac{\rho}{\rho_o}q_o -
        \rho \int_{\rho_o}^{\rho} \frac{c(r)}{r}\,dr
        & \rho\ge\rho_o
      \end{array}
    \right.
    \\
    L_2^-(\rho;u_o)
    & = &
    \left\{
      \begin{array}{l@{\quad\mbox{ if }}l}
        \frac{\rho}{\rho_o}q_o +
        \sqrt{\frac{\rho}{\rho_o} (\rho-\rho_o) 
          \left(p(\rho)-p(\rho_o)\right)}
        & \rho\ge\rho_o,\\
        \frac{\rho}{\rho_o}q_o +
        \rho \int_{\rho_o}^{\rho} \frac{c(r)}{r}\,dr
        & \rho\le\rho_o,
      \end{array}
    \right.
  \end{array}
\end{equation}

%
%

\begin{proofof}{Proposition~\ref{prop:Local}}
  A direct application of the Implicit Function Theorem allows to
  define neighborhoods $\mathcal{A}_l, \mathcal{A}_r$ of $\bar a_l,
  \bar a_r$ and $\mathcal{U}_l, \mathcal{U}_r$ of $\bar u_l, \bar u_r$
  such that for all $a_l \in \mathcal{A}_l$, $a_r \in \mathcal{A}_r$
  and $u_l \in \mathcal{U}_l$, there exists a unique $u_r \in
  \mathcal{U}_r$ such that $\Phi(a_l,u_l;a_r,u_r) = 0$ if and only if
  $u_r = T(u_l;a_l,a_r)$.

  The proof of {1.--6.}~is extremely similar to various results
  already appeared in the literature. More precisely, {1.--5.}~are
  shown as in~\cite[Theorem~3.3]{ColomboGaravello3},
  \cite[Theorem~3.2]{ColomboHertySachers},
  \cite[Theorem~2.3]{ColomboMarcellini}.
  The latter estimate follows
  from~\cite[formula~(3.2)]{ColomboHertySachers}
  and~\textbf{($\mathbf{\Phi}$\ref{it:p0})}.
\end{proofof}

\begin{proofof}{Proposition~\ref{prop:allP}}
  Concerning~\textbf{($\mathbf{\Phi}$\ref{it:p0})}, in the different
  cases, simple computations show that
  \begin{displaymath}
    D_{u_r} \Phi (a_l,u_l;a_r,u_r)
    =
    \left\{
      \begin{array}{l@{\qquad\mbox{ in case }}c}
        -{a_r}^2 \, \lambda_1(u_r) \, \lambda_2(u_r)
        &
        \textbf{(L)}
        \\
        a_r \,p'(\rho_r)
        &
        \textbf{(p)}
        \\
        -a_r \, \lambda_1(u_r) \, \lambda_2(u_r)
        &
        \textbf{(P)}
        \\
        -a_r \, \lambda_1(u_r) \, \lambda_2(u_r)
        &
        \textbf{(S)}
      \end{array}
    \right.
  \end{displaymath}
  which does not vanish in $\pint{A}_0$.

  In case~\textbf{(S)},
  consider~\textbf{($\mathbf{\Phi}$\ref{it:p1})}. The condition $\Phi
  (a, \rho_l, q_l; a, \rho_r, q_r) = 0$ implies $q_l = q_r$ and so
  $\rho_l = \rho_r$; hence~\textbf{($\mathbf{\Phi}$\ref{it:p1})}
  holds.  Concerning~\textbf{($\mathbf{\Phi}$\ref{it:p2})} and assume
  that $\Phi, (a_l, \rho_l, q_l; a_r, \rho_r, q_r) = 0$, which is
  equivalent to
  \begin{displaymath}
    \left\{
      \begin{array}{l}
        a_l q_l = a_r q_r,\\
        \rho_r = R (a_r; \rho_l, q_l).
      \end{array}
    \right.
  \end{displaymath}
  Since $a Q(a) = a_l q_l$ for every $a$ and the equation for $R$
  depends on $a_l^2 q_l^2$, then we deduce that, if $\rho_r = R (a_r;
  \rho_l, q_l)$, then $\rho_l = R (a_l; \rho_r, q_r) = R (a_l; \rho_r,
  - q_r)$.  Thus $\Phi (a_l, \rho_l, q_l; a_r, \rho_r, q_r) = 0$ is
  equivalent to $\Phi (a_r, \rho_r, - q_r; a_l, \rho_l, - q_l) = 0$
  and~\textbf{($\mathbf{\Phi}$\ref{it:p2})} holds.  The fact that $a
  Q(a) = a_l q_l$ for every $a$ and the equation for $R$ depends on
  $a_l^2 q_l^2$ implies that~\textbf{($\mathbf{\Phi}$\ref{it:p3})}
  also holds. Property~\textbf{($\mathbf{\Phi}$\ref{it:p4})} follows
  from~2.~in Lemma~\ref{le:5.1}.

  The other cases are immediate.
\end{proofof}

\begin{proofof}{Theorem~\ref{thm:Lsol}}
  Assume that there exists an~\textbf{(L)}-solution $u$ attaining
  values in $A_0$. Then
  \begin{displaymath}
    (\rho(t,0-), L_1(\rho(t,0-); \bar u_l)) \in A_0, \quad
    (\rho(t,0+), L_2^-(\rho(t,0+); \bar u_r)) \in A_0
  \end{displaymath}
  and $l' \leq \rho(t,0-)$, $l'' \ge \rho(t,0-)$. Moreover,
  by~\cite[Lemma~1]{ColomboGaravello1}, we have
  \begin{eqnarray*}
    & & a_l P \left(l', L_1(l'; \bar u_l)\right) - 
    a_r P\left (g \left(\frac{a_l}{a_r} L_1(l'; \bar u_l) \right),
      \frac{a_l}{a_r} L_1(l'; \bar u_l)\right)
    \\
    & \leq &
    a_l P\left(u(t,0-)\right) - a_r P\left(u(t,0+)\right)  =   0
    \,,\qquad \mbox{ and }
    \\
    & & a_l P \left(l'', L_1(l''; \bar u_l)\right) - a_r P\left (g
      \left(\frac{a_l}{a_r} L_1(l''; \bar u_l) \right),
      \frac{a_l}{a_r} L_1(l''; \bar u_l)\right)
    \\
    & \ge & a_l P\left(u(t,0-)\right) - a_r P\left(u(t,0+)\right) = 0
  \end{eqnarray*}
  Hence~(\ref{eq:L-sol_CNS}) is a necessary condition.

  Assume now that~(\ref{eq:L-sol_CNS}) holds. Since $l'\leq l''$, then
  \begin{displaymath}
    (\rho, L_1(\rho; \bar u_l)) \in A_0, \quad
    \left(g \left(\frac{a_l}{a_r} L_1(\rho; \bar u_l) \right),
      \frac{a_l}{a_r} L_1(\rho; \bar u_l)\right) \in A_0
  \end{displaymath}
  for every $\rho \in [l', l'']$. Moreover the function
  \begin{displaymath}
    \rho \longmapsto
    a_l P \left(\rho, L_1(\rho; \bar u_l)\right) - 
    a_r P\left (g \left(\frac{a_l}{a_r}
        L_1(\rho; \bar u_l) \right),
      \frac{a_l}{a_r} L_1(\rho; \bar u_l)\right),
  \end{displaymath}
  defined in $[l',l'']$, is increasing,
  by~\cite[Lemma~1]{ColomboGaravello1}, and continuous. Hence there
  exists $\tilde \rho \in [l',l'']$ such that
  \begin{displaymath}
    a_l P \left(\tilde \rho, L_1(\tilde \rho; \bar u_l)\right) - 
    a_r P\left (g \left(\frac{a_l}{a_r}
        L_1(\tilde \rho; \bar u_l) \right),
      \frac{a_l}{a_r} L_1(\tilde \rho; \bar u_l)\right) = 0.
  \end{displaymath}
  Therefore the traces of an~\textbf{(L)}-solution at $J$ are given by
  \begin{displaymath}
    (\tilde \rho, L_1(\tilde \rho; \bar u_l)) \quad \textrm{ and }
    \quad
    \left(g \left(\frac{a_l}{a_r}
        L_1(\tilde \rho; \bar u_l) \right), 
      \frac{a_l}{a_r} L_1(\tilde \rho; \bar u_l) \right).
  \end{displaymath}
  This concludes the existence proof.

  Assume that the Riemann Problem~(\ref{eq:RP}) admits two
  different~\textbf{(L)}-solu\-tions, attaining values in $A_0$,
  denoted with $(\rho', q')$ and $(\rho'', q'')$. Note that
  \begin{displaymath}
    \begin{array}{rcl@{\qquad}rcl}
      q'(t,0-) & = & L_1(\rho'(t,0-); \bar u_l) \,,
      &
      q''(t,0-) & = & L_1(\rho''(t,0-); \bar u_l) \,,
      \\
      q'(t,0+) & = & L_2^-(\rho'(t,0+); \bar u_r) \,,
      &
      q''(t,0+) & = & L_2^-(\rho''(t,0+); \bar u_r) \,.
    \end{array}
  \end{displaymath}
  Without loss of generality, we may assume that $\rho'(t,0+) <
  \rho''(t,0+)$.  By~\cite[Lemma~1]{ColomboGaravello1}, we deduce that
  $q'(t,0+) < q''(t,0+)$ and
  \begin{displaymath}
    a_r P(\rho'(t,0+), q'(t,0+)) < a_r P(\rho''(t,0+), q''(t,0+)).
  \end{displaymath}
  By definition of~\textbf{(L)}-solution, the previous inequality
  becomes
  \begin{displaymath}
    a_l P(\rho'(t,0-), q'(t,0-)) < a_l P(\rho''(t,0-), q''(t,0-)),
  \end{displaymath}
  and, by~\cite[Lemma~1]{ColomboGaravello1}, we deduce that
  $\rho'(t,0-) < \rho''(t,0-)$ and so $q'(t,0-) > q''(t,0-)$, which is
  in contradiction with $q'(t,0+) < q''(t,0+)$.
\end{proofof}

\begin{proofof}{Theorem~\ref{th:p-sol}}
  Assume that a~\textbf{(p)}-solution $(\rho, q)$ in $A_0$ exists.
  Then
  \begin{displaymath}
    (\rho(t,0-), L_1(\rho(t,0-); \bar u_l)) \in A_0, \quad
    (\rho(t,0+), L_2^-(\rho(t,0+); \bar u_r)) \in A_0
  \end{displaymath}
  and so $l' \leq \rho(t,0-) = \rho(t,0+) \leq l''$. Moreover,
  by~\cite[Lemma~1]{ColomboGaravello1}, we have
  \begin{eqnarray*}
    & &
    a_l \, L_1(l'; \bar u_l) -
    a_r \, L_2^-(l'; \bar u_r) \\
    & \ge &
    a_l \, L_1(\rho(t,0-); \bar u_l) -
    a_r \, L_2^-(\rho(t,0+); \bar u_r)  =   0 \,.
    \\
    & & a_l \, L_1(l''; \bar u_l) -
    a_r \, L_2^-(l''; \bar u_r) \\
    & \leq & a_l \, L_1(\rho(t,0-); \bar u_l) - a_r L_2^-(\rho(t,0+);
    \bar u_r) = 0 \,.
  \end{eqnarray*}
  Hence~(\ref{eq:p-sol_CNS}) is a necessary condition.

  Assume now that~(\ref{eq:p-sol_CNS}) holds.  First, we claim that
  $l'\leq l''$. Suppose, by contradiction, that $l' > l''$. We have
  two different possibilities: either $\fd{l} (\bar u_l) < \fd{r}
  (\bar u_r)$ or $\fu{r} (\bar u_r) < \fu{l} (\bar u_l)$.  Consider
  only the first case, the second one being similar.  We easily deduce
  that $l' = \fd{r} (\bar u_r)$ and $l'' = \fd{l} (\bar u_l)$ and so
  \begin{equation}
    \label{eq:ineq_L1_L2}
    L_1 (l'; \bar u_l) <
    L_2^- (l'; \bar u_r)
    \quad \mbox{ and } \quad
    L_1 (l''; \bar u_l) >
    L_2^- (l''; \bar u_r).
  \end{equation}
  Moreover, (\ref{eq:p-sol_CNS}) implies that
  \begin{displaymath}
    \frac{L_2^- (l''; \bar u_r)}{L_1 (l''; \bar u_l)}
    \le \frac{a_l}{a_r} \le
    \frac{L_2^- (l'; \bar u_r)}{L_1 (l'; \bar u_l)},
  \end{displaymath}
  which is in contradiction with~(\ref{eq:ineq_L1_L2}). Since $l'\leq
  l''$, then
  \begin{displaymath}
    (\rho, L_1(\rho; \bar u_l)) \in A_0, \quad
    (\rho, L_2^-(\rho; \bar u_r)) \in A_0
  \end{displaymath}
  for every $\rho \in [l', l'']$. Moreover the function
  \begin{displaymath}
    \rho \longmapsto a_l \, L_1(\rho; \bar u_l) -
    a_r \, L_2^-(\rho; \bar u_r),
  \end{displaymath}
  defined in $[l',l'']$, is decreasing,
  by~\cite[Lemma~1]{ColomboGaravello1}, and continuous. Hence there
  exists $\tilde \rho \in [l',l'']$ such that
  \begin{displaymath}
    a_l \, L_1(\tilde \rho; \bar u_l) -
    a_r \, L_2^-(\tilde \rho; \bar u_r) = 0.
  \end{displaymath}
  Therefore the traces of a~\textbf{(p)}-solution at $J$ are given by
  \begin{displaymath}
    (\tilde \rho, L_1(\tilde \rho; \bar u_l)) \quad \textrm{ and }
    \quad
    (\tilde \rho, L_2^-(\tilde \rho; \bar u_r)).
  \end{displaymath}

  Assume now that the Riemann Problem~(\ref{eq:RP}) admits two
  different~\textbf{(p)}-solutions $(\rho', q')$ and $(\rho'', q'')$
  attaining values in $A_0$. Thus, we deduce that
  \begin{displaymath}
    \rho' (t, 0+) = \rho' (t, 0-) \ne \rho'' (t, 0-) = \rho'' (t, 0+) 
  \end{displaymath}
  for a.e. $t>0$. Without loss of generality, we suppose that
  \begin{displaymath}
    \rho' (t, 0+) = \rho' (t, 0-) < \rho'' (t, 0-) = \rho'' (t, 0+) 
  \end{displaymath}
  for a.e. $t>0$ and so, by using~\cite[Lemma~1]{ColomboGaravello1},
  \begin{displaymath}
    \begin{array}{rclcl}
      a_r \, L_2^- (\rho' (t,0+); \bar u_r)
      & = &
      a_l \, L_1 (\rho' (t,0-); \bar u_l)
      & > &
      a_l \, L_1 (\rho'' (t,0-); \bar u_l)\\
      & = &
      a_r \, L_2^- (\rho'' (t,0+); \bar u_r)
      & > &
      a_r \, L_2^- (\rho' (t,0+); \bar u_r),
    \end{array}
  \end{displaymath}
  for a.e. $t>0$, which is a contradiction.  This concludes the proof.
\end{proofof}

\begin{proofof}{Corollary~\ref{cor:Nop}}
  Consider only the case $\fu{r} (\bar u_r) < \fu{l} (\bar u_l)$, the
  other one being similar. Assume by contradiction that $(\rho, q)$ is
  a~\textbf{(p)}-solution to~(\ref{eq:RP}) attaining values in
  $A_0$. Therefore, we deduce that
  \begin{displaymath}
    \begin{array}{rcl}
      q(t,0-) & = & L_1(\rho (t,0-);\bar u_l),
      \\
      q(t,0+) & = & L_2^-(\rho (t,0+);\bar u_r),
      \\
      p \left(\rho (t,0-) \right) & = & p \left(\rho (t,0+) \right),
    \end{array}
    \qquad \mbox{for a.e. } t > 0 \,,
  \end{displaymath}
  hence $\rho (t,0-) = \rho (t,0+)$ for a.e.~$t > 0$.  Moreover, $\rho
  (t,0-) \geq \fu{l} (\bar u_l)$ and $\rho (t,0+) \leq \fu{r} (\bar
  u_r)$, which gives a contradiction and proves non existence.
\end{proofof}

\begin{proofof}{Corollary~\ref{cor:NoP}}
  Consider only the case $\fu{r} (\bar u_r) < \fu{l} (\bar u_l)$, the
  other one being similar. Assume by contradiction that $(\rho, q)$ is
  a~\textbf{(P)}-solution to~(\ref{eq:RP}), attaining values in
  $A_0$. Therefore we deduce that
  \begin{displaymath}
    \begin{array}{rcl}
      q(t,0-) & = & L_1(\rho (t,0-);\bar u_l),\\
      q(t,0+) & = & L_2^-(\rho (t,0+);\bar u_r),\\
      P(\rho (t,0-), q(t,0-)) & = & P(\rho (t,0+), q(t,0+)),
    \end{array}
  \end{displaymath}
  for a.e. $t>0$. Moreover we have
  \begin{eqnarray*}
    P(\rho (t,0-), q(t,0-)) & \ge &
    P(\fu{l} (\bar \rho_l, 0),
    L_1(\fu{l} (\bar \rho_l, 0); \bar u_l))\\
    & > & P(\fu{r} (\bar \rho_r, 0),
    L_2^-(\fu{r} (\bar \rho_r, 0); \bar u_r))\\
    & \ge & P(\rho (t,0+), q(t,0+))
  \end{eqnarray*}
  and so we get a contradiction, proving not existence.
\end{proofof}

\begin{proofof}{Lemma~\ref{le:5.1}}
  The first item is trivial; so we consider only the second one.

  If $q_l = 0$, then $\frac{d}{da}R(a)=0$ and $R(a)$ is constantly
  equal to $\rho_l$. Assume therefore $q_l \ne 0$. The sign of the
  derivative of $R(a)$ is given by the sign of
  \begin{displaymath}
    a^2 p'(R(a)) R^2(a) - a_l^2 q_l^2.
  \end{displaymath}
  At the initial point $a_l$, this term is $a_l^2 \left( p'(\rho_l)
    \rho_l^2 - q_l^2\right)$, which is strictly positive if and only
  if $(\rho_l, q_l) \in A_0$. Thus $R(a)$ is increasing in a
  neighborhood of $a_l$. If $R(a)$ is increasing, then we deduce that
  \begin{displaymath}
    a^2 p'(R(a)) R^2(a) - a_l^2 q_l^2 \ge a_l^2
    \modulo{p'(\rho_l) \rho_l^2 - q_l^2}
  \end{displaymath}
  for every $a \ge a_l$, provided $R(a)$ exists.  Using this estimate
  and the comparison theorem for ODE, we conclude that $R(a)$ exists
  and is increasing for every $a > a_l$.  Finally, for $a > a_l$
  sufficiently big, the derivative $R'(a)$ can be bounded by $K /
  a^3$, where $K$ is a constant depending on the initial
  conditions. Thus we deduce that $R(a)$ is bounded.
\end{proofof}

\begin{lemma}
  \label{le:derivative_R}
  Fix $(\rho_l, q_l) \in A_0$ and denote with $R(a; \rho_l, q_l)$ the
  solution to the ODE in~(\ref{eq:sysRQ_explicit}).  Define
  \begin{displaymath}
    z_1(a) := \frac{\partial}{\partial \rho}
    R(a; \rho, q_l)_{\vert \rho = \rho_l} \quad \textrm{ and } \quad
    z_2(a) := \frac{\partial}{\partial q}
    R(a; \rho_l, q)_{\vert q = q_l}.
  \end{displaymath}
  Then $z_1$ and $z_2$ satisfy the following system
  \begin{equation}
    \label{eq:z-1-2}
    \left\{
      \begin{array}{l}
        \frac{d}{da}z_1(a) = \frac{\partial}{\partial \rho}\, 
        g\left(a, R(a; \rho_l, q_l), q_l)\right)\,z_1(a)\vspace{.2cm}\\
        \frac{d}{da}z_2(a) = \frac{\partial}{\partial \rho} 
        g\left(a, R(a; \rho_l, q_l), q_l)\right) z_2(a)\!
        + \! \frac{\partial}{\partial q}
        g\left(a, R(a; \rho_l, q_l), q_l)\right)\vspace{.2cm}\\
        z_1(a_l) = 1\vspace{.2cm}\\
        z_2(a_l) = 0,
      \end{array}
    \right.
  \end{equation}
  where $\displaystyle g \left(a, \rho, q\right) =
  \frac{\rho}{a}\,\frac{a_l^2 q^2}{a^2 p'(\rho) \rho^2 - a_l^2 q^2}$.
\end{lemma}

The proof consists in the classical derivation of the solution to an
ODE with respect to a parameter, hence we omit it.

\begin{proofof}{Theorem~\ref{th:S-sol}}
  Consider the following curve on $A_0$
  \begin{equation}
    \label{eq:curve_psi}
    \begin{array}{rccc}
      \psi : & \left[ \fu{l}(\bar \rho_l, \bar q_l),
        \fd{l}(\bar \rho_l, \bar q_l) \right] & \longrightarrow
      & A_0\vspace{.2cm}\\
      & s & \longmapsto & (\psi_1(s), \psi_2(s)),
    \end{array}
  \end{equation}
  where $\psi_1(s) = R (a_r; s, L_1(s; \bar \rho_l, \bar q_l))$ and
  $\psi_2(s) = \frac{a_l}{a_r} L_1(s; \bar \rho_l, \bar q_l))$.  Since
  the point $(s, L_1(s; \bar \rho_l, \bar q_l))$ belongs to $A_0$ for
  every $s$ in the domain of $\psi$, then we easily deduce that the
  image of $\psi$ is contained in $A_0$, by Lemma~\ref{le:5.1}.
  Clearly the second component of $\psi$ is decreasing with respect to
  $s$.

  We claim that $\psi_1$ is increasing with respect to $s$.  In the
  following we use the same notation of Lemma~\ref{le:derivative_R}.
  The derivative of $\psi_1$ is
  \begin{eqnarray*}
    \psi_1'(s) & = & \frac{d}{ds}\, R (a_r; s, L_1(s; \bar \rho_l, \bar q_l))\\
    & = & \frac{\partial}{\partial s}\,
    R (a_r; s, L_1(s; \bar \rho_l, \bar q_l)) +
    \frac{\partial}{\partial q}\, R (a_r; s, L_1(s; \bar \rho_l, \bar q_l))
    \frac{\partial}{\partial s}\,L_1(s; \bar \rho_l, \bar q_l)\\
    & = & z_1(a_r; s, L_1(s; \bar \rho_l, \bar q_l)) + 
    z_2(a_r; s, L_1(s; \bar \rho_l, \bar q_l)) 
    \frac{\partial}{\partial s}\,L_1(s; \bar \rho_l, \bar q_l).
  \end{eqnarray*}
  By~(\ref{eq:z-1-2}), we have that $z_1(a_r; s, L_1(s; \bar \rho_l,
  \bar q_l)) > 0$ for every $s$ in the domain of the curve $\psi$.
  Let us consider some different cases.
  \begin{enumerate}
  \item $L_1(s, \bar \rho_l, \bar q_l) > 0$. In this case the
    derivative
    \begin{displaymath}
      \frac{\partial}{\partial q} 
      g(a_r, R(a_r; s, L_1(s, \bar \rho_l, \bar q_l)),
      L_1(s, \bar \rho_l, \bar q_l))
    \end{displaymath}
    is strictly positive. For $a \ge a_l$, define the function
    \begin{displaymath}
      \beta (a) = z_1(a; s, L_1(s; \bar \rho_l, \bar q_l)) + 
      z_2(a; s, L_1(s; \bar \rho_l, \bar q_l)) 
      \frac{\partial}{\partial s}\,L_1(s; \bar \rho_l, \bar q_l).
    \end{displaymath}
    Easy computations show that
    \begin{displaymath}
      \left\{
        \begin{array}{rcl}
          \frac{d}{da} \beta (a) & = & 
          \frac{\partial}{\partial \rho} 
          g(a, R(a; s, L_1(s, \bar \rho_l, \bar q_l)),
          L_1(s, \bar \rho_l, \bar q_l)) \beta(a) \vspace{.2cm}\\
          & & + \frac{\partial}{\partial q} 
          g(a, R(a; s, L_1(s, \bar \rho_l, \bar q_l)),
          L_1(s, \bar \rho_l, \bar q_l)),\vspace{.2cm}\\
          \beta (a_l) & = & 1.
        \end{array}
      \right.
    \end{displaymath}
    Define
    \begin{displaymath}
      \bar a = \inf \left \{ a \ge a_l : \beta(a) = 0 \right \}.
    \end{displaymath}
    Assume by contradiction that $\bar a < +\infty$.  In this case we
    deduce that
    \begin{displaymath}
      \frac{d}{da} \beta (\bar a) = \frac{\partial}{\partial q} 
      g(\bar a, R(\bar a; s, L_1(s, \bar \rho_l, \bar q_l)),
      L_1(s, \bar \rho_l, \bar q_l)) > 0
    \end{displaymath}
    by assumptions and this is not possible. Hence $\beta(a) > 0$ for
    every $a \ge a_l$. In particular $\beta (a_r) > 0$ and so
    $\psi_1'(s) > 0$.

  \item $L_1(s, \bar \rho_l, \bar q_l) = 0$. In this case the
    derivative
    \begin{displaymath}
      \frac{\partial}{\partial q} 
      g(a_r, R(a_r; s, L_1(s, \bar \rho_l, \bar q_l)),
      L_1(s, \bar \rho_l, \bar q_l))
    \end{displaymath}
    vanishes and so, by~(\ref{eq:z-1-2}), $z_2(a_r; s, L_1(s; \bar
    \rho_l, \bar q_l)) = 0$.  Hence $\psi_1'(s) > 0$.

  \item $L_1(s, \bar \rho_l, \bar q_l) < 0$. In this case the
    derivative
    \begin{displaymath}
      \frac{\partial}{\partial q} 
      g(a_r, R(a_r; s, L_1(s, \bar \rho_l, \bar q_l)),
      L_1(s, \bar \rho_l, \bar q_l))
    \end{displaymath}
    is strictly negative and so, by~(\ref{eq:z-1-2}), $z_2(a_r; s,
    L_1(s; \bar \rho_l, \bar q_l)) < 0$.  Hence $\psi_1'(s) > 0$.
  \end{enumerate}

  \noindent By the previous considerations, a necessary and sufficient
  condition for the existence and uniqueness of
  an~\textbf{(S)}-solution $(\rho, q)$, attaining values in $A_0$, is
  that the image of $\psi$ intersects in a unique point the image of
  the curve
  \begin{equation}
    \label{eq:curve2}
    \begin{array}{ccc}
      [\fd{r}(\bar u_r), \fu{r}(\bar u_r)] & \longrightarrow & A_0\\
      \rho & \longmapsto & \left( \rho, L_2^-(\rho; \bar u_r) \right).
    \end{array}
  \end{equation}
  The image of the curve~(\ref{eq:curve2}) divides the set $A_0$ in
  two parts and this permits to conclude.
\end{proofof}

\begin{proofof}{Corollary~\ref{cor:S_noexistence}}
  Assume by contradiction that $(\rho, q)(t,x)$ is
  an~\textbf{(S)}-solution to~(\ref{eq:RP}), attaining values in
  $A_0$.  Therefore we deduce that
  \begin{displaymath}
    \begin{array}{rcl}
      q(t,0-) & = & L_1(\rho (t,0-);\bar u_l),
      \\
      q(t,0+) & = & L_2^-(\rho (t,0+);\bar u_r),
      \\
      R \left(\alpha_r; \rho (t,0-), q(t,0-) \right) & = & \rho (t,0+),
    \end{array}
  \end{displaymath}
  for a.e. $t>0$. By Lemma~\ref{le:5.1}, we deduce that $R(\alpha_r;
  \rho (t,0-), q(t,0-)) > \rho (t,0-)$.  Moreover, by hypotheses, we
  have that
  \begin{displaymath}
    \rho (t,0-) \ge \fu{l} (\bar u_l) >
    \fu{r} (\bar\rho_r, \bar q_r) \ge \rho (t,0+)
  \end{displaymath}
  and so we obtain a contradiction.
\end{proofof}

{\small

  {\bibliographystyle{abbrv}
    
    \bibliography{psj}}}

\end{document}